\documentclass{amsart}
\usepackage{tabu}
\usepackage[margin=95pt]{geometry}
\usepackage{enumitem}
\allowdisplaybreaks

\usepackage{amssymb}
\usepackage{amsmath}
\usepackage{amscd}
\usepackage{amsbsy}
\usepackage{comment}
\usepackage{enumitem}
\usepackage[matrix,arrow]{xy}
\usepackage{hyperref}
\def\mod#1{{\ifmmode\text{\rm\ (mod~$#1$)}
\else\discretionary{}{}{\hbox{ }}\rm(mod~$#1$)\fi}}

\begin {document}

\newtheorem{thm}{Theorem}
\newtheorem{lem}{Lemma}[section]
\newtheorem{prop}[lem]{Proposition}

\newtheorem{cor}[lem]{Corollary}

\theoremstyle{definition}

\theoremstyle{remark}

\title[Arithmetic progressions]
{Arithmetic progressions in sumsets of geometric progressions}

\author[Michael Bennett]{Michael A. Bennett}
\address{Department of Mathematics, University of British Columbia, Vancouver, B.C., V6T 1Z2 Canada}
\email{bennett@math.ubc.ca}

\thanks{The author is supported by NSERC}

\date{\today}


\begin {abstract}
If $a$ and $b$ are integers with $b>a>1$, we completely characterize ``long'' arithmetic progressions in the sumsets of the
geometric progressions $1, a, a^2, a^3, \ldots$ and $1, b, b^2, b^3, \ldots$. Our proofs utilize recent applications of bounds for  linear forms in logarithms to $S$-unit equations, and consequences of the modularity of  Frey-Hellegouarch curves, together with elementary arguments. 
\end {abstract}
\maketitle

%------------------------------
\section{Introduction}
%------------------------------

Let $a$ and $b$ be integers with $b>a>1$, and define the set $S_{a,b}$ to be the sumset of the geometric progressions
$$
1, a, a^2, a^3, \ldots \; \mbox{ and } \; 1, b, b^2, b^3, \ldots,
$$
i.e.
$$
S_{a,b} = \{ n \in \mathbb{N} \; : \; n = a^x+b^y, \; \mbox{ for } x, y \in \mathbb{Z}, \; x, y \geq 0 \}.
$$
The main goal  of this paper is to study the arithmetic progressions in the sets $S_{a,b}$. These sets are sufficiently thin that, were they random, our expectation for a fixed pair $(a,b)$  would be that $S_{a,b}$ should contain at most finitely many three-term arithmetic progressions, and indeed this is the case ``most'' of the time. By way of example, $S_{2,7}$ contains precisely  $22$ such progressions, while, in contrast, $S_{2,9}$ contains infinitely many  and $S_{2,8}$ even contains infinitely many four-term arithmetic progressions. It is likely that $S_{a,b}$ contains infinitely many 
three-term arithmetic progressions precisely when either $a=2$ and $b=2^k+1$, for $k$ a positive integer, or when $a$ and $b$ are multiplicatively  dependent, i.e. there exist positive integers $c$ and $d$ such that $a^c=b^d$. We will discuss this further in Section \ref{Sec9}.

In this paper, we will focus our attention on the problem of characterizing longer arithmetic progressions in the sets $S_{a,b}$. 
Our starting point  is recent work of  Chen, Huang and Zhang \cite{ChHuZh} who showed, via a completely elementary argument,  that the longest arithmetic progression in the set 
$S_{2,3}$
has length $6$. Regarding such maximal progressions, we prove the following broad generalization of this, characterizing all $5$-term progressions in the sets $S_{a,b}$.

\begin{thm} \label{thm-main}
Let $a$ and $b$ be integers with $b>a>1$. Suppose that there exist positive integers $N$ and $D,$ and nonnegative integers $x_i, y_i$ such that
$$
N+iD = a^{x_i}+b^{y_i}, \; \; \mbox{ for } \; i \in \{ 0, 1, 2, 3, 4 \}.
$$
Then we have either
$$
(a,b,N,D)=(2,2^k+1,2^k+1,2^k)  \; \mbox{ or } \;  (3,4 \cdot 3^{k-1}+1,3^{k-1}+1,2 \cdot 3^{k-1}) \; \mbox{ for } \; k \in \mathbb{N},
$$
or
$$
\begin{array}{c}
(a,b,N,D) \in \left\{ (2,3,5,2),  (2,3,7,6), (2,3,9,8),  (2,3,17,24),  \right. \\
\left. (2,3,41,24), (2,5,5,8),  (2,9,17,24),  (2,9,41,24), (3,4,7,6) \right\}.
\end{array}
$$
\end{thm}

\noindent An almost immediate corollary of this is the following.
\begin{cor} \label{cor-main}
If $a$ and $b$ are integers with $b>a>1$,  then the only $6$-term arithmetic progressions in $S_{a,b}$ are 
$$
3,5,7,9,11,13 \; \; \mbox{ and } \; \; 17, 41, 65, 89, 113, 137,
$$
in $S_{2,3}$ and in both $S_{2,3}$ and $S_{2,9}$, respectively.
\end{cor}

Our arguments, in contrast to those of  \cite{ChHuZh}, rely, implicitly or explicitly, upon various results from Diophantine approximation, including bounds for linear forms in logarithms, both complex and $p$-adic, and upon various Diophantine consequences of the modularity of Frey-Hellegouarch curves. The outline of this paper is the following.
In Section \ref{Sec2}, we present the various results we will need for explicitly solving certain polynomial-exponential and $3, 4$ and $5$-term $S$-unit equations. Sections \ref{Sec3}--\ref{Sec7} are devoted to dealing with $5$-term arithmetic progressions in $S_{a,b}$ with fixed numbers of ``large'' terms, proving Theorem \ref{thm-main}. In Section \ref{Sec8}, we prove Corollary \ref{cor-main}. Finally, in Section \ref{Sec9}, we discuss various families of pairs $(a,b)$ for which the $S_{a,b}$ are known to have $3$ or $4$-term progressions.

%--------------------------------------------------------------------------------
\section{Notation and preliminary results} \label{Sec2}
%--------------------------------------------------------------------------------

Suppose that $a$ and $b$ are integers with $b>a \geq 2$ and that we have a $5$-term arithmetic progression in integers of the shape $a^x+b^y$,  i.e. that there exist positive integers $N$ and $D,$ and nonnegative integers $x_i, y_i$ such that
$$
N+iD = a^{x_i}+b^{y_i}, \; \; \mbox{ for } 0 \leq i \leq 4.
$$
We will write
$$
n = \left[ \frac{\log (N+4D)}{\log a} \right] \; \mbox{ and } \; m = \left[ \frac{\log (N+4D)}{\log b} \right],
$$
where by $[x]$ we mean the greatest integer $\leq x$.
Let us call a term $N+iD$ in our arithmetic progression {\it $a$-large} if there is some representation of $N+iD = a^{x_i}+b^{y_i}$ with $x_i=n$, {\it $b$-large} if there is a representation of $N+iD = a^{x_i}+b^{y_i}$ with $y_i=m$, and {\it large}, if  $N+iD$ is either $a$-large or $b$-large.

Note for future use  that, if $k$ and $j$ are nonnegative integers,
\begin{equation} \label{better}
a^{n-k}+b^{m-j} <   \left( \frac{1}{a^k}+\frac{1}{b^j} \right) (N+4D).
\end{equation}
In particular, from $b > a \geq 2$, applying this with $k=j=1$, $N+4D$ is necessarily large.

We will have need of the following technical lemma, based primarily upon  Theorem 1.2 of an old paper of the author \cite{Be-BLMS}. It is perhaps worth noting that these results rely fundamentally upon Frey-Hellegouarch type arguments applied to ternary equations of the shape $ax^n+by^n=cz^2$.

\begin{lem} \label{tech-lem}
If $b>2$ is an integer, and $x, y, \alpha$ and $\beta$ are nonnegative integers such that
$$
b^x-b^y=2^\alpha 3^\beta,
$$
then
\begin{equation} \label{fammy}
b=2^\alpha 3^\beta+1, \; x=1, \; y=0,
\end{equation}
\begin{equation} \label{fammy2}
b=3, \; x=y+1, \; \beta=y, \;  \alpha=1,
\end{equation}
\begin{equation} \label{fammy3}
b=3, \; x=y+2, \; \beta=y, \;  \alpha=3,
\end{equation}
\begin{equation} \label{fammy4}
b=9, \; x=y+1, \; \beta=2y, \;  \alpha=3,
\end{equation}
\begin{equation} \label{fammy5}
b=4, \; x=y+1, \; \alpha=2y, \;  \beta=1,
\end{equation}
or $y=0$ and 
\begin{equation} \label{spec}
(b,x,\alpha,\beta) \in \{ (2,2,0,1), (5,2,3,1), (7,2,4,1), (17,1,5,2) \}.
\end{equation}

\end{lem}
\begin{proof}
Necessarily, we have $x > y \geq 0$. If $y=0$, then Theorem 1.2 and Corollary 1.4 of \cite{Be-BLMS} imply that $x=1$ and so we have (\ref{fammy}), or that $(b,x,\alpha,\beta)$ are as in (\ref{spec}). 
If $y > 0$, then $b=2^{\alpha_1} 3^{\beta_1}$ for integers $\alpha_1, \beta_1$ with
$$
\alpha_1 y \leq \alpha \; \mbox{ and } \; \beta_1 y \leq  \beta.
$$
If $\alpha_1 \beta_1 > 0$, then we have $\alpha_1 y = \alpha$ and $\beta_1 y =  \beta$, contradicting $b > 2$. If $\alpha_1=0$, we have $\beta_1 y = \beta$ and so
$$
3^{\beta_1 (x-y)} - 1 = 2^\alpha,
$$
whence either $\beta_1 (x-y)=\alpha=1$ or $\beta_1 (x-y) =2, \;\alpha=3$. These correspond to (\ref{fammy2}),  (\ref{fammy3}), and (\ref{fammy4}).
If $\beta_1=0$, then $\alpha_1y=\alpha$ and
$$
2^{\alpha_1 (x-y)}-1 = 3^\beta.
$$
Thus $\alpha_1(x-y)=2$ and $\beta=1$. Since $b>2$, we necessarily have $\alpha_1=2$ and $x=y+1$,  corresponding to (\ref{fammy5}).
\end{proof}

We will also have need of a variety of special results on $3$, $4$ and $5$-term $S$-unit equations. The first, for $5$-term equations where the set $S= \{ 2, 3 \}$, is a special case of Theorem 8 of \cite{BajBen}. The proof of this result relies upon careful application of bounds for linear forms in complex and $p$-adic logarithms, together with combinatorial arguments.

\begin{prop}[Bajpai and B.] \label{BajBen}
Suppose that $\alpha_i, \beta_j$ are nonnegative integers and that
\begin{equation} \label{biggie}
2^{\alpha_1} 3^{\beta_1} \pm 2^{\alpha_2} 3^{\beta_2} \pm 2^{\alpha_3} 3^{\beta_3} \pm 2^{\alpha_4} 3^{\beta_4} \pm 2^{\alpha_5} 3^{\beta_5} =0,
\end{equation}
where
$$
\gcd ( 2^{\alpha_1} 3^{\beta_1}, 2^{\alpha_2} 3^{\beta_2}, 2^{\alpha_3} 3^{\beta_3}, 2^{\alpha_4} 3^{\beta_4}, 2^{\alpha_5} 3^{\beta_5} ) =1
$$
and (\ref{biggie}) has no {\it vanishing subsums}, i.e. for every $1 \leq  i < j \leq 5$ we have $2^{\alpha_i} 3^{\beta_i} \neq 2^{\alpha_j} 3^{\beta_j}$. Then
$$
\max_{1 \leq i \leq 5} \{ 2^{\alpha_i} 3^{\beta_i} \} \leq 3^{12}, \;
\max_{1 \leq i \leq 5} \{ \alpha_i \} \leq 19 \; \mbox{ and } \;
\max_{1 \leq i \leq 5} \{ \beta_i \} \leq 12.
$$
\end{prop}

For $S$ containing slightly larger primes, we will 
appeal to Theorem 6.3 of  de Weger \cite{DeWe} (for $3$-term equations) and both Theorems 1 and 2 of Deze and Tijdeman \cite{DeTi} (for $4$-term equations). 

\begin{prop}[de Weger]  \label{DW}
The Diophantine equation
$$
x+y=z
$$
has precisely $545$ solutions in positive integers $x, y$ and $z$ with $\gcd (x,y)=1$, $x \leq y$ and where the greatest prime factor $P(x \cdot y \cdot z) \leq 13$, each satisfying
$$
\begin{array}{c}
\mbox{ord}_2 (x \cdot y \cdot z) \leq  15, \; \; \mbox{ord}_3 (x \cdot y \cdot z) \leq  10, \; \; \mbox{ord}_5 (x \cdot y \cdot z) \leq  7, \; \; \\
\mbox{ord}_7 (x \cdot y \cdot z) \leq  6, \; \; \mbox{ord}_{11} (x \cdot y \cdot z) \leq  5, \; \; \mbox{ord}_{13} (x \cdot y \cdot z) \leq  4. \; \; \\
\end{array}
$$
Here, by $\mbox{ord}_p (m)$, we mean the largest power of a prime $p$ which divides a nonzero integer $m$.
\end{prop}

\begin{prop}[Deze and Tijdeman] \label{DT}
Let $p$ and $q$ be distinct primes with $\max \{ p, q \} < 200$. Then if we have either
\begin{equation} \label{DT1}
p^xq^y \pm p^z \pm q^w \pm 1=0
\end{equation}
or 
\begin{equation} \label{DT2}
p^x \pm q^y \pm p^z \pm q^w = 0,
\end{equation}
then 
$$
\max \{ p^x, q^y, p^z, q^w \} \leq 2^{15}.
$$
\end{prop}

%-----------------------------------------------------------
\section{Three equal exponents} \label{Sec3}
%-----------------------------------------------------------

The basic idea behind the elementary arguments in \cite{ChHuZh} is that having sufficiently many terms of the shape $2^{x_i}+3^{y_i}$ in arithmetic progression guarantees that several of the exponents either coincide or differ by a small amount. We generalize this argument to $S_{a,b}$; 
our first result towards proving Theorem \ref{thm-main} is the following.

\begin{prop} \label{equality}
Let $a$ and $b$ be integers with $b>a>1$. Suppose that there exist positive integers $N$ and $D,$ and nonnegative integers $x_i, y_i$ such that
$$
N+iD = a^{x_i}+b^{y_i}, \; \; \mbox{ for } \; i \in \{ 0, 1, 2, 3, 4 \}.
$$
Suppose further that here exist integers  $i, j$ and $k$ with $0 \leq i < j < k \leq 4$ and either 
\begin{equation} \label{x}
x_i=x_j=x_k,
\end{equation}
or
\begin{equation} \label{y}
y_i=y_j=y_k.
\end{equation}
Then either
$$
(a,b,N,D)=(2,3,5,2), (2,3,7,6), (2,3,9,8), (3,4,7,6),
$$
\begin{equation} \label{family}
(a,b,N,D)=(2,2^{x_0}+1,2^{x_0}+1,2^{x_0})
\end{equation}
or
\begin{equation} \label{family2}
(a,b,N,D)=(3,4 \cdot 3^{x_0}+1,3^{x_0}+1,2 \cdot 3^{x_0}).
\end{equation}
\end{prop}

Before we begin, we would like to emphasize how our arguments in the remainder of this section and in Sections \ref{Sec4} -- \ref{Sec7} will proceed.
Various linear relations between $3$ or more of the terms $N+iD$, $i \in \{ 0, 1, 2, 3, 4 \}$, lead immediately to $S$-unit equations in the primes dividing $ab$. Repeated application of inequalities like (\ref{better})
enable us to derive small absolute upper bounds for  certain differences of exponents $|x_i-x_j|$ and $|y_k-y_l|$ (equations (\ref{x}) and (\ref{y}) are, in a certain sense, best case scenarios of this phenomenon), which in turn reduce the number of terms in the $S$-unit equations. Further, straightforward inequalities like 
$$
4(N+D)>N+4D \; \mbox{ and } \; 2(N+2D) > N+4D, 
$$
combined with these upper bounds, lead to absolute bounds upon $b$ (and hence $a$), whereby we can treat the remaining equations with a combination of Lemma \ref{tech-lem} and Propositions \ref{BajBen}, \ref{DW} and \ref{DT}. Our claim is that a careful combinatorial analysis leads to this desired conclusion. Regrettably for  the author and undoubtedly for the reader, this analysis takes the form of an extremely unpleasant, but perhaps unavoidable, case-by-case argument. Here, we should stress, that this approach apparently does require the assumption of $5$ terms in an arithmetic progression, and fails to lead to a similar conclusion, at least in full generality, in characterizing $4$-term progressions in $S_{a,b}$.

\begin{proof}[Proof of Proposition \ref{equality}]
Suppose first that (\ref{x}) holds with $x_i=x_j=x_k=s$.
The identity
$$
 (k-j) (N+iD)+(j-i)(N+k D) =(k-i)(N+jD)
 $$
 thus implies that we have one of 
$$
b^{y_i}+b^{y_k} = 2 b^{y_j}, \;  b^{y_i}+2b^{y_k} = 3 b^{y_j},   \; b^{y_k}+2b^{y_i} = 3 b^{y_j}, b^{y_k}+3b^{y_i} = 4 b^{y_j}  \mbox{ or } b^{y_i}+3b^{y_k} = 4 b^{y_j},
$$
where $y_i < y_j < y_k$. Dividing by $b^{y_i}$, we contradict $b \geq 3$ unless $b^{y_k}+3b^{y_i} = 4 b^{y_j}$, corresponding to $i=0, j=1, k=4$. 
In this  case $b=3$, $y_1=y_0+1$, $y_4=y_0+2$ and
necessarily $a=2$, so that
$$
N=2^s + 3^{y_0}, \; N+D = 2^s + 3^{y_0+1}, \; N+4D = 2^s + 3^{y_0+2},
$$
whence $D=2 \cdot 3^{y_0}$ and we have
\begin{equation} \label{Baj}
2^{x_3}+3^{y_3} = 2^{x_2}+3^{y_2} + 2 \cdot 3^{y_0}.
\end{equation}
We first need to consider the possibility of vanishing subsums in this equation. These correspond to one of $x_2=x_3$, $y_2=y_3$,
$x_3=y_2=0$, $y_3=x_2=0$, or $x_3=1, y_0=0$.

If $x_2 = x_3$, necessarily, since then $y_2<y_3$, we have $y_2=y_0$ and hence $y_3=y_2+1=y_0+1=y_1$. We thus have
$$
D = 2 \cdot 3^{y_0}= 2^{x_2-1}-2^{s-1},
$$
so that
$$
2^{x_2-2}-2^{s-2} = 3^{y_0}
$$
and thus $s=2$ and either $x_2=3$ and $y_0=0$, or $x_2=4$ and $y_0=1$. These correspond to
$$
(a,b,N,D)=(2,3,5,2) \; \mbox{ and } \; (2,3,7,6),
$$
respectively. Similarly, if $y_2=y_3$, then 
$$
2^{x_3-1} = 2^{x_2-1}+ 3^{y_0}
$$
and so $x_2=1$ and either $y_0=0$ and $x_3=2$, or $y_0=1$ and $x_3=3$. From
$$
D = 2 \cdot 3^{y_0}= (N+2D)-(N+D)=3^{y_2}-2^s-1,
$$
we thus have
$$
3^{y_2} - 2 \cdot 3^{y_0} =2^s+1.
$$
Appealing to Proposition \ref{DT}, we find that  $s=1$, $y_0=1$ and $y_2=2$, whence
$$
N+2D=2+3^2 = N+D,
$$
a contradiction.

Next, suppose that $x_3=y_2=0$ in equation (\ref{Baj}). Then
$$
3^{y_3} = 2^{x_2}+ 2 \cdot 3^{y_0}
$$
and so $x_2=0$ which implies that $N+2D=2$, a contradiction. Similarly, $y_3=x_2=0$ leads to 
$$
2^{x_3}=3^{y_2} + 2 \cdot 3^{y_0},
$$
again a contradiction via parity. Finally, $x_3=1$ and $y_0=0$ yields 
$$
3^{y_3} = 2^{x_2}+3^{y_2},
$$
whereby $y_2=0$ and $x_2 \in \{ 1, 3 \}$. But this gives that either $N+2D =3$ (contradicting $D =2 \cdot 3^{y_0}=2$) or $N+2D=9$ (so that $(a,b,N,D)=(2,3,5,2)$).

We may thus suppose that (\ref{Baj}) has no vanishing subsums, whereby we may appeal to Proposition \ref{BajBen} to conclude that $x_3 \leq 19$ and $y_3 \leq 12$; a short computation verifies that there are precisely $2$ tuples  $(x_3,y_3,x_2,y_2,y_0)$ satisfying  equation (\ref{Baj}), with no vanishing subsums and the additional constraint that $2^{x_2}+3^{y_2}-5 \cdot 3^{y_0}$ is a power of $2$ :
$$
(x_3,y_3,x_2,y_2,y_0) =(3,0,2,1,0) \mbox{ and } (1,2,3,0,0).
$$
These correspond to $(a,b,N,D)=(2,3,3,2)$ and $(2,3,5,2)$, respectively.

%(all with $\max \{ x_3,y_3,x_2,y_2,y_0 \} \leq 10$). 
%$$
%(x_2,y_2) \in \{ (0,0), (0,1), (0,2), (0,3), (0,4), (1,0),  (1,2), (2,1), (2,3),  (3,1), (3,5), (5,3),  (9,3) \}.
%$$

Similarly, if (\ref{y}) holds with $y_i=y_j=y_k=t$, 
 then we have one of 
$$
a^{x_i}+a^{x_k} = 2 a^{x_j}, \;  a^{x_i}+2a^{x_k} = 3 a^{x_j},   \; a^{x_k}+2a^{x_i} = 3 a^{x_j}, a^{x_k}+3a^{x_i} = 4 a^{x_j}  \mbox{ or } a^{x_i}+3a^{x_k} = 4 a^{x_j},
$$
where $x_i < x_j < x_k$. Dividing by $a^{x_i}$, we contradict $a \geq 2$ unless we have one of  
$$
a^{x_k}+2a^{x_i} = 3 a^{x_j},  \; i=0, j=1, k=3,
$$
or
$$
a^{x_k}+2a^{x_i} = 3 a^{x_j},  \; i=1, j=2, k=4,
$$
or
$$
a^{x_k}+3a^{x_i} = 4 a^{x_j}, \; i=0, j=1, k=4.
$$
In the first two cases, $a=2$, while in the third $a=3$. We thus have 
\begin{equation} \label{case1}
N = 2^{x_0}+b^t, \; N+D = 2^{x_0+1}+b^t, \; N+3D=2^{x_0+2}+b^t,
\end{equation}
\begin{equation} \label{case2}
N+D = 2^{x_1}+b^t, \; N+2D = 2^{x_1+1}+b^t, \; N+4D=2^{x_1+2}+b^t,
\end{equation}
or
\begin{equation} \label{case3}
N=3^{x_0} +b^t, \; N+D = 3^{x_0+1} + b^t, \; N+4D = 3^{x_0+2} + b^t.
\end{equation}

We begin by treating case (\ref{case1}), where $D=2^{x_0}$, 
\begin{equation} \label{e4}
N+4D= 5 \cdot 2^{x_0}+b^t = 2^{x_4}+b^{y_4}
\end{equation}
and
\begin{equation} \label{e2}
N+2D = 3 \cdot 2^{x_0}+b^t = 2^{x_2}+b^{y_2}.
\end{equation}
If $y_4=y_2$, subtracting these last two equations yields
$$
2^{x_4} = 2^{x_0+1} + 2^{x_2},
$$
whence $x_2=x_0+1$ and $x_4=x_0+2$. But then
$$
b^{y_2}-b^t = 2^{x_0}
$$
and so, from Lemma \ref{tech-lem},  $t=0$ and either
$$
y_2=1, \; b=2^{x_0}+1 \; \mbox{ or } \; y_2=2, \; b=x_0=3,
$$
corresponding to (\ref{family}) and $(a,b,N,D)=(2,3,9,8)$, respectively. We may thus suppose that $y_4 \neq y_2$ and, from (\ref{e4}) and (\ref{e2}), $y_4 \neq t$ and $y_2 \neq t$.

Since $N+4D$ is large,  we necessarily have $x_4=n$ or $y_4=m$. If $t=m$, it follows that $y_2, y_4 < m$  and $x_4=n$, and so, from
$$
2(N+2D) = 2^{x_2+1} + 2 \cdot b^{y_2}  = 2^{x_0}+b^m + 2^{n}+b^{y_4} > 2^n + b \cdot b^{y_2},
$$
we have that  $x_2 =n$, whereby
$$
b^{y_4}-b^{y_2} = 2^{x_0+1}.
$$
Again appealing to Lemma \ref{tech-lem},
$y_2=0$, and either $y_4=1$ and $b=2^{x_0+1}+1$, or $y_4=2$, $b=3$ and $x_0=2$.
In the first case, we have
$$
2^{x_0} = D = (N+2D)-(N+D) = 2^n+1-2^{x_0+1}-b^m \equiv 0 \mod{2^{x_0+1}},
$$
a contradiction. 
If, on the other hand, $y_4=2$, $b=3$ and $x_0=2$, 
$$
2^{n}  =3^m +11,
$$
a contradiction modulo $8$.

We may thus assume that $t < m$, whereby
$$
b^t \leq \frac{1}{b} b^m < \frac{1}{b} (N+4D)
$$
and so, from (\ref{e4}), 
$$
5 \cdot 2^{x_0} > \left( 1 - \frac{1}{b} \right) (N+4D) > \left( 1 - \frac{1}{b} \right) 2^n.
$$
Since $b \geq 3$ and $N+4D > 5 \cdot 2^{x_0}$, we thus conclude that $n=x_0+2$. If $x_4=n$, it follows from (\ref{e4}), that 
$$
b^{y_4}-b^t=2^{x_0},
$$
whence, again from Lemma \ref{tech-lem},  either
$$
y_4=1, \; t=0, \; b=2^{x_0}+1,
$$
or
$$
b^{y_4}=9, \; t=0, \; x_0=3.
$$
The first of these corresponds to (\ref{family}),  while the second has $(a,b,N,D)=(2,3,9,8)$ or $(2,9,9,8)$. If, on the other hand, $x_4 < n$, then $y_4=m$ and so $y_2 < m$. From $2(N+2D) > N+4D$, 
$$
 2^{x_2+1} > \left( 1 - \frac{2}{b} \right) 2^{x_0+2},
 $$
 and so, since $b \geq 3$, necessarily $x_2 \geq x_0$, whence $x_2 \in \{ x_0, x_0+1, x_0+2 \}$. Equation (\ref{e2}) thus implies that
 $$
 \left| b^{y_2}-b^t \right| = 2^\kappa, \; \; \kappa \in \{ x_0, x_0+1 \}.
 $$
 Once again, we have that $t=0$ or that $y_2=0$; no new arithmetic progressions  accrue.
 
In case (\ref{case2}), we have $D=2^{x_1}$, whence 
$$
N=b^t=2^{x_0}+b^{y_0}.
$$
It follows from Lemma \ref{tech-lem} that $y_0=0$ and either $b=3$, $t=2$ and $x_0=3$, or $b=2^{x_0}+1$, $t=1$. In the first case, 
$$
N+3D = 9 + 3 \cdot 2^{x_1} = 2^{x_3}+3^{y_3},
$$
so that $y_3=0$, whereby
$$
8 + 3 \cdot 2^{x_1} = 2^{x_3},
$$
and so $x_1=3$ and $x_3=5$, corresponding to $(a,b,N,D)=(2,3,9,8)$.
If, on the other hand, $b=2^{x_0}+1$ and $t=1$, then necessarily $y_3=0$ and we have
$$
N+3D = 2^{x_0}+1 + 3 \cdot 2^{x_1} = 2^{x_3}+1,
$$
whence
$$
2^{x_0} + 3 \cdot 2^{x_1} = 2^{x_3}.
$$
Thus $x_0=x_1$ and we are in case (\ref{family}).

Finally, suppose we are in situation (\ref{case3}), whence $a=3$. Then $D=2 \cdot 3^{x_0}$ and we have
\begin{equation} \label{eqy1}
N+2D=5 \cdot 3^{x_0} + b^t = 3^{x_2}+b^{y_2}
\end{equation}
and
\begin{equation} \label{eqy2}
N+3D=7 \cdot 3^{x_0} + b^t = 3^{x_3}+b^{y_3}.
\end{equation}
If $y_2=y_3$, then $3^{x_3}-3^{x_2} = 2 \cdot 3^{x_0}$, so that $x_2=x_0$ and $3^{x_3-x_2}-1 = 2$, i.e. $x_3=x_2+1$.
Thus
$$
b^{y_3}-b^t=4 \cdot 3^{x_0}.
$$
Lemma \ref{tech-lem} then implies that either $b=4, y_3=2, t=1, x_0=1$ (corresponding to $(a,b,N,D)=(3,4,7,6)$), or that $t=0, y_3=1, b=4 \cdot 3^{x_0}+1$ (corresponding to (\ref{family2})).

Similarly, if $x_2=x_3$, then
$$
b^{y_3}-b^{y_2}= 2 \cdot 3^{x_0}.
$$
Once again from Lemma \ref{tech-lem} and the fact that $b \geq 4$, we have $y_2=0, \; y_3=1$ and $b=2 \cdot 3^{x_0}+1$, contradicting  equation (\ref{eqy1}) modulo $8$.

We may thus suppose that $y_2 \neq y_3$ and that $x_2 \neq x_3$, so that, in particular, at least one of $y_2, y_3$ is distinct from $m$, say $y_i < m$, $i \in \{ 2, 3 \}$. We have
$$
N+iD < \left( \frac{1}{3^{m-x_i}} + \frac{1}{b^{m-y_i}} \right) (N+4D)
$$
and so either $x_i=n$, or we have $i=2$, $x_2=n-1$, $y_2=m-1$ and $b \in \{ 4, 5 \}$. 
In these latter cases, necessarily
$$
N+2D=3^{n-1}+b^{m-1} \; \mbox{ and } \; N+3D= 3^n+b^{y_3}.
$$
If $t<m$, then, since $N+4D$ is large, $x_0=n-2$, and thus, from $N+(N+4D)=2(N+2D)$,
$$
3^{n-2}  +3^{n} + 2 b^t = 2 (3^{n-1}+b^{m-1}),
$$
i.e.
$$
b^{m-1}-b^t = 2 \cdot 3^{n-2}. 
$$
Lemma \ref{tech-lem} and the fact that $b \in \{ 4, 5 \}$ imply that $m=2$, a contradiction. If $t=m$, 
$$
b^m-b^{m-1}+3^{x_0+2}-3^{n-1}= 4 \cdot 3^{x_0}
$$
so that
$$
(b-1) \cdot b^{m-1}= 3^{n-1} -5 \cdot 3^{x_0},
$$
a contradiction modulo $b$, for $b \in \{ 4, 5 \}$. 

We may therefore assume that  $x_i=n$ and that also $y_{5-i}=m$. If $t=m$, we thus have
$$
3^{x_0+2} - 3^{x_{5-i}} = (i-1) \cdot 2 \cdot 3^{x_0},
$$
i.e.
$$
3^{x_{5-i}} = (11-2i) \cdot 3^{x_0},
$$
an immediate contradiction. We have $t<m$ and hence again 
$x_0=n-2$, whence
\begin{equation} \label{final}
(N+4D) -(N+iD) = b^t-b^{y_i}=(4-i)D=(4-i) \cdot 2 \cdot 3^{n-2}.
\end{equation}
From Lemma \ref{tech-lem}, we have one of
$$
i=2, \; b=4, \; t=2, \; y_2=1, \; n=3,
$$
$$
i=2, \; y_2=0, \; t=1, \; b=4 \cdot 3^{n-2}+1,
$$
or
$$
i=3, \; y_3=0, \; t=1, \; b=2 \cdot 3^{n-2}+1.
$$
In the first case, $N=19$ and $D=6$, so that $N+3D=37$, a contradiction.
In the second, 
$$
N+4D = 3^n+b < \frac{13}{4} b \leq \frac{13}{16} b^2 < b^m < N+3D,
$$
while, in the third we have
$$
N+4D = 3^n+b < \frac{11}{2} b,
$$
and
$$
N+2D > b^m \geq b^2,
$$
so that $b \in \{ 4, 5 \}$, contradicting $b=2 \cdot 3^{n-2}+1$.
This completes the proof of Proposition \ref{equality}.
\end{proof}

%-----------------------------------------------
\section{After Proposition \ref{equality} : four large terms} \label{Sec4}
%-----------------------------------------------

With Proposition \ref{equality} in hand, to complete the proof of Theorem \ref{thm-main}, it remains to show that the only $5$-term arithmetic progressions in any of the $S_{a,b}$, where both (\ref{x}) and (\ref{y}) fail to hold, are those given by
$$
(a,b,N,D) \in \{ (2,3,17,24), (2,3,41,24),  (2,9,17,24), (2,9,41,24) \}.
$$
Let us therefore suppose that we have
$$
N+iD = a^{x_i}+b^{y_i}, \; \; \mbox{ for } \; i \in \{ 0, 1, 2, 3, 4 \},
$$
but that (\ref{x}) and (\ref{y}) are not satisfied for any $0 \leq i < j < k \leq 4$. In particular, at least one of these terms, say  $N+\kappa D$ is necessarily  not large, whence
\begin{equation} \label{starty}
N+\kappa D \leq a^{n-1}+b^{m-1}.
\end{equation}
Also, necessarily $m \geq 2$ (else $y_i \in \{ 0, 1 \}$ for all $i$ and hence we have three equal values of the exponents $y_i$).

We begin, in this section, by considering  the case where four of $N, N+D, N+2D, N+3D$ and $N+4D$ large,
 whereby, from our assumptions, precisely two of them are $a$-large and the other two are $b$-large, with no term simultaneously $a$-large and $b$-large. In particular, there exist indices $0 \leq i < j \leq 4$ with
 $x_i=x_j =n$ and  $y_i< y_j < m$, and indices $0 \leq k < l \leq 4$ with $y_k=y_l =m$, $x_k< x_l < n$, and $\{i,j,k,l,\kappa \} = \{ 0,1,2,3,4 \}$. It follows that 
\begin{equation} \label{deee}
D = \frac{1}{l-k} \left( a^{x_l}-a^{x_k} \right) = \frac{1}{j-i} \left( b^{y_j}-b^{y_i} \right),
\end{equation} 
so that, from (\ref{starty}), 
$$
N+4D \leq a^{n-1} + b^{m-1} + (4-\kappa) \frac{1}{j-i} \left( b^{y_j}-b^{y_i}  \right) < a^{n-1} + \left( 1+ \frac{4-\kappa}{(j-i) \cdot b^{m-1-y_j}} \right) b^{m-1}.
$$
We thus have
$$
N+4D < \left( \frac{1}{a} + \frac{1}{b} + \frac{4-\kappa}{(j-i) \cdot b^{m-y_j}} \right) (N+4D)
$$
and hence
\begin{equation} \label{goop1}
\frac{1}{a} + \frac{1}{b} + \frac{4-\kappa}{(j-i) \cdot b^{m-y_j}} > 1.
\end{equation}
Arguing similarly, we find that 
\begin{equation} \label{goop2}
\frac{1}{a} + \frac{1}{b} + \frac{4-\kappa}{(l-k) \cdot a^{n-x_l}} > 1.
\end{equation}
In particular, from (\ref{goop1}), 
$$
2 \leq a \leq 5 \; \mbox{ and } \; b < \frac{5a}{a-1}.
$$
From (\ref{goop1}) and (\ref{goop2}), we further have 
\begin{equation} \label{poop}
m-2 \leq y_j \leq m-1 \; \mbox{ and } \; n-4 \leq x_l \leq n-1,
\end{equation}
and, if $a \in \{ 4, 5 \}$, necessarily
\begin{equation} \label{poop2}
y_j = m-1, \; x_l = n-1,\; j-i=l-k=1, \; \kappa=0 \; \mbox{ and } \; N=a^{n-1}+b^{m-1}.
\end{equation}
If $a \in \{ 4, 5 \}$, it follows that 
$$
\frac{1}{j} \left( a^n-a^{n-1} \right) = a^{n-1}-a^{x_k},
$$
a contradiction unless $a=j=4$ (so that $(i,j,k,l)=(3,4,1,2)$) and $x_1=n-2$. Since, in this case,
$$
\frac{1}{2} \left( b^m-b^{m-1} \right) = b^{m-1}-b^{y_i},
$$
we contradict $b \geq 5$. We thus may assume that 
$$
a \in \{ 2, 3 \} \; \mbox{ and } \; b < \frac{5a}{a-1}.
$$ 
For these cases, if
$$
(a,b,i,j,k,l,\kappa) \not\in \{ (2, 5, 2, 3, 1, 4, 0), 
(2, 7, 2, 3, 1, 4, 0), (3, 5, 3, 4, 0, 2, 1) \},
$$
from inequalities (\ref{goop1}) and (\ref{goop2}), equation (\ref{deee}) in each case reduces to an equation of the shape (\ref{DT2}), whereby we can apply  Proposition \ref{DT} to conclude, in each case, that 
$$
\max \left\{ \frac{1}{l-k} a^{x_l}, \frac{1}{j-i} b^{y_j} \right\} \leq 2^{15}.
$$
A short computation reveals only the solutions to
$$
p^x-p^y=q^z-q^w
$$
with $(p,q,x,y,z,w)$ one of
$$
\begin{array}{l}
(2,3,2,1,1,0), (2,3,3,1,2,1), (2,3,4,3,2,0), (2,3,5,3,3,1), (2,3,8,4,5,1), (2,5,3,2,1,0),  \\
(2,5,5,3,2,0), (2,5,7,2,3,0), (2,5,7,3,3,1), (2,7,3,1,1,0), (2,7,6,4,2,0), (3,5,3,1,2,0),  \\
\end{array}
$$
or $(3,7,2,1,1,0)$. None of these correspond to new progressions $(a,b,N,D)$.

If, on the other hand,
$$
(a,b,i,j,k,l,\kappa) \in \{ (2, 5, 2, 3, 1, 4, 0), 
(2, 7, 2, 3, 1, 4, 0), (3, 5, 3, 4, 0, 2, 1) \},
$$
we necessarily have, from (\ref{goop1}), that $y_j=m-1$. If, further, $N+\kappa D = a^{n-1}+b^{m-1}$, then 
$$
\frac{1}{j-\kappa} a^{n}-a^{n-1} =  \frac{1}{l-k} \left( a^{x_l}-a^{x_k} \right) \leq  \frac{1}{l-k} \left( a^{n-1}-a^{n-2} \right),
$$
in each case a contradiction. 
Thus $N+\kappa D \neq a^{n-1}+b^{m-1}$, whence 
$$
N+\kappa D \leq \max \{ a^{n-2}+b^{m-1}, a^{n-1}+b^{m-2} \}.
$$
In particular,
$$
\frac{1}{a} + \frac{1}{b^2} + \frac{4-\kappa}{(j-i) \cdot b^{m-y_j}} > 1.
$$
This is a contradiction for $(a,b,i,j,k,l,\kappa)=(2, 7, 2, 3, 1, 4, 0)$ or  $(3, 5, 3, 4, 0, 2, 1)$. We may thus suppose that $(a,b,i,j,k,l,\kappa)=(2, 5, 2, 3, 1, 4, 0)$. In this case, $y_3=m-1$ and 
$$
D = \frac{1}{3} \left( 2^{x_4}-2^{x_1} \right) = 5^{m-1}-5^{y_2},
$$
so that, in particular, $x_4 \geq n-3$. But then, from
$$
D=(N+4D)-(N+3D)= 2^{x_4}+5^m-2^n-5^{m-1},
$$
we have
$$
5^{m-1}-5^{y_2} = 2^{x_4}+5^m-2^n-5^{m-1}
$$
whereby
$$
3 \cdot 5^{m-1}+5^{y_2} = 2^n-2^{x_4}.
$$
Since 
$$
3 \cdot 5^{m-1}-3 \cdot 5^{y_2} = 2^{x_4}-2^{x_1},
$$
we have
$$
5^{y_2} = 2^{n-2}-2^{x_4-1}+2^{x_1-2},
$$
and so $x_1=2$. Since $2^{x_4} \equiv 2^{x_1} \mod{3}$, necessarily $x_4 \geq x_1+2 \geq 4$ and so, modulo $8$, $y_2$ is even. We may thus appeal to Theorem 2 of Szalay \cite{Sz} to conclude that $y_2=2$, $x_4=4$ and $n=7$. But then $D=4$ and so $m=2$, whence $2^7 < N+4D < 5^3$, a contradiction.

%----------------------------------------------------------
\section{Three large terms} \label{Sec5}
%----------------------------------------------------------

Suppose next that precisely three of $N+iD$, $i \in \{ 0, 1, 2, 3, 4 \}$ are large. From our previous work, we may assume that either
\begin{equation} \label{Case1}
\mbox{there exist } 0 \leq i < j \leq 4 \mbox{ and } k \not\in \{ i, j \}  \mbox{ with } x_i=x_j=n, \; y_k=m, 
\end{equation}
or
\begin{equation} \label{Case2}
\mbox{there exist } 0 \leq i < j \leq 4 \mbox{ and } k \not\in \{ i, j \}  \mbox{ with } y_i=y_j=m, \; x_k=n.
\end{equation}
Since two of the terms $N+iD$, $i \in \{ 0, 1, 2, 3, 4 \}$ are not large, in all cases we have inequality (\ref{starty}) with $\kappa=1$. 
In particular,
$$
N+D \leq a^{n-1}+b^{m-1} < \left( \frac{1}{a} + \frac{1}{b} \right) (N+4D)
$$
and hence, from $4(N+D)>N+4D$, $\frac{1}{a} + \frac{1}{b} > \frac{1}{4}$,
i.e.
\begin{equation} \label{special2}
a \in \{ 2, 3, 4 \}, \; a=5, \; 6 \leq b \leq 19, \; \;  a=6, \; 7 \leq b \leq 11 \; \mbox{ or } \; a=7, \; 8 \leq b \leq 9.
\end{equation}

Suppose first that $N+4D$ is both $a$-large and $b$-large, so that, from
$$
N+D \leq a^{n-1}+b^{m-1} < \frac{1}{a} (N+4D),
$$
necessarily $a \in \{ 2, 3 \}$. If either of $N+2D$ or $N+3D$ is not large, then
$$
N+2D \leq a^{n-1}+b^{m-1} < \frac{1}{a} (N+4D),
$$
contradicting $a \geq 2$. We thus have either
\begin{equation} \label{dcase1}
D = \frac{1}{2} \left( a^n-a^{x_2} \right) = b^m-b^{y_3}
\end{equation}
or
\begin{equation} \label{dcase2}
D =a^n-a^{x_3}= \frac{1}{2} \left(  b^m-b^{y_2} \right).
\end{equation}
In case (\ref{dcase1}), if $a^n \leq 2 b^m$, then 
$$
N=a^n+b^m - 4 \left( b^m-b^{y_3} \right) \leq a^n -3b^m+4 b^{m-1}
\leq (4-b) \, b^{m-1}
$$
and hence $b=3$, so that $a=2$. Equation (\ref{dcase1}) thus becomes
$$
2^{n-1}-2^{x_2-1} = 3^m-3^{y_3}
$$
and so, from Proposition \ref{DT}, we have $n \leq 16$. A short computation gives that 
$$
(n,x_2,m,y_3) \in \{  (3,2,1,0), (4,2,2,1), (5,4,2,0), (6,4,3,1), (9,5,5,1) \},
$$
and so, since $N>1$, 
$$
(a,b,N,D) = (2,3,3,2) \; \mbox{ or } \;   (2,3,9,8).
$$

If we have (\ref{dcase1}) and $a^n > 2 b^m$, 
$$
N=a^n+b^m - 2  \left( a^n-a^{x_2} \right) < 2 a^{x_2} - \frac{1}{2} a^n
$$
and so either $a=2$ and $x_2 \geq n-2$, or $a=3$ and $x_2=n-1$.
In the first case, (\ref{dcase1}) becomes either
$$
b^m-b^{y_3} = 2^{n-2}
$$
or
$$
b^m-b^{y_3} = 3 \cdot 2^{n-3}.
$$
In the second, we have
$$
b^m-b^{y_3} = 3^{n-1}.
$$
Applying Lemma \ref{tech-lem}, the first of these equations leads after a little work to (\ref{family}), while the second and third, in each case, contradicts $N>1$.

If, on the other hand, we are in case  (\ref{dcase2}), if $b^m \geq 2 a^n$, then
$$
N=a^n+b^m - 2 \left( b^m-b^{y_3} \right) \leq 2b^{y_3} - \frac{1}{2} b^m \leq \frac{1}{2} (4-b) \, b^{m-1},
$$
whence $(a,b)=(2,3)$. But then
$$
N+4D > b^m \geq 2a^n =2^{n+1}> N+4D.
$$
We may thus suppose that $b^m < 2 a^n$, whence 
$$
N=a^n+b^m -   4 \left( a^n-a^{x_3} \right) < 4 a^{x_3} - a^n.
$$
It follows that, again, either $a=2$ and $x_3 \geq n-2$, or $a=3$ and $x_3=n-1$, whereby we have
$$
b^m-b^{y_2}  = 2^n,
$$
$$
b^m-b^{y_2} = 3 \cdot 2^{n-1} 
$$
or
$$
b^m-b^{y_2} = 4 \cdot 3^{n-1}.
$$
From Lemma \ref{tech-lem}, we in each case contradict $N>1$.

We may thus suppose that $N+4D$  fails to be both $a$-large and $b$-large (whereby the same is true for all $N+iD$). It follows from (\ref{Case1}) and (\ref{Case2})  that either
\begin{equation} \label{hope1}
(j-i) D =  b^{y_j} - b^{y_i}, \; y_j \leq m-1
\end{equation}
or that
\begin{equation} \label{hope2}
(j-i) D =  a^{x_j} - a^{x_i}, \; x_j \leq n-1.
\end{equation}
In either case, we will assume that $N+\kappa D$ is not large, so that we may write
$$
N+\kappa D = a^{n-\delta_1}+b^{m-\delta_2}, \; \delta_i \geq 1,
$$
and have
$$
\{ i,j,k,\kappa \} = \{ 1, 2, 3, 4 \}.
$$
In particular, we have
$$
N+\kappa D < \frac{1}{a^{\delta_1}} (N+kD) +  \frac{1}{b^{\delta_2}} (N+iD) 
$$
This gives a contradiction unless 
\begin{equation} \label{awesome}
\frac{k}{a^{\delta_1}} + \frac{i}{b^{\delta_2}} > \kappa.
\end{equation}

\subsection{Two $a$-large terms}

In case (\ref{hope1}), it follows that 
\begin{equation} \label{foosh!}
D < \frac{1}{(j-i)b^{m-y_j}} (N+kD)
\end{equation}
whence
\begin{equation} \label{ook!}
N+\kappa D > \left( \frac{(j-i)b^{m-y_j}+\kappa-k}{(j-i)b^{m-y_j}+4-k}. \right) (N+4D).
\end{equation}
Thus
$$
\left( \frac{1}{a^{\delta_1}} +\frac{1}{b^{\delta_2}} \right) (N+4D) > 
a^{n-\delta_1}+b^{m-\delta_2} = N+\kappa D >\left( \frac{(j-i)b^{m-y_j}+\kappa-k}{(j-i)b^{m-y_j}+4-k} \right) (N+4D)
$$
and so
\begin{equation} \label{spooky}
 \frac{1}{a^{\delta_1}} +\frac{1}{b^{\delta_2}} >  \frac{(j-i)b^{m-y_j}+\kappa-k}{(j-i)b^{m-y_j}+4-k}.
 \end{equation}
 This inequality provides an upper bound upon $b$ (and hence $a$) and, unless we have 
 \begin{equation} \label{real-special}
 (i,j,k,\kappa) =(2,3,4,1), \; y_j=m-1, \;  \mbox{ and } \; (a,b)=(2,3), 
 \end{equation}
 upon $\min \{ \delta_1, \delta_2 \}$. In particular, a short computation allows us to sharpen (\ref{special2}) to conclude that, necessarily
 \begin{equation} \label{special3}
 (a,b) \in \{ (2,3), (2,4), (2,5), (2,6), (2,7), (3,4), (3,5), (4,5) \}.
 \end{equation}
  
 The identities
 $$
 (\kappa-j) (N+iD)+(j-i)(N+\kappa D) +(i-\kappa)(N+jD)=0
 $$
 and
 $$
 (\kappa-j) (N+kD)+(j-k)(N+\kappa D) +(k-\kappa)(N+jD)=0
 $$
 lead  to the equations
  \begin{equation} \label{gaew1}
\left(  (j-i) a^{\delta_1}   + i-j \right) a^{n-\delta_1} = (\kappa -j) b^{y_i} + (j-i) b^{m-\delta_2} + (i-\kappa ) b^{y_j}
 \end{equation} 
 and
   \begin{equation} \label{gaew2}
\left(  (\kappa -j) b^{\delta_2}   + j-\kappa \right) b^{m-\delta_2} +(k-\kappa) b^{y_j} = (j-\kappa ) a^{x_k} + (k-j) a^{n-\delta_1} + (\kappa-k ) a^n,
 \end{equation} 
 respectively. 
 The first of these describes a $5$-term $S$-unit equation with $S = \{ 2, 3 \}$, in case 
 $$
 (a,b) \in \{ (2,3), (2,4),  (2,6),  (3,4) \}.
 $$
 Appeal to Proposition \ref{BajBen}  bounds $n$ and $m$, and a short calculation uncovers no new progressions.

From (\ref{special3}), we may thus suppose that
  \begin{equation} \label{special4}
 (a,b) \in \{ (2,5), (2,7),  (3,5), (4,5) \}.
 \end{equation}
 From (\ref{spooky}), it follows that, in all cases, $y_j=m-1$. Further, if $(a,b)=(2,5)$, then we have either
 $$
 (i,j,k,\kappa,\delta_1,\delta_2)=(3,4,2,1,1,1)
 $$
 or
 $$
 (i,j,k,\kappa)=(2,3,4,1) \; \mbox{ and either } \; \delta_1=1 \;
 \mbox{ or } \; (\delta_1,\delta_2)=(2,1).
 $$
 For the remaining three pairs $(a,b)=(2,7), (3,5)$ and $(4,5)$, necessarily
 $$
 (i,j,k,\kappa,\delta_1,\delta_2)=(2,3,4,1,1,1).
 $$

 Further,
 \begin{equation} \label{gaew3}
 D = \frac{1}{j-i} \left( b^{y_j}-b^{y_i} \right) = \frac{1}{|k-j|} \left| a^{x_k} + b^{m}-a^n-b^{y_j} \right|.
 \end{equation}
 
If $\delta_1=1$ and 
$$
(a,b) \in \{ (2,5), (2,7), (3,5), (4,5) \},
$$
equation (\ref{gaew1}) becomes
\begin{equation} \label{mum}
(a-1) \cdot a^{n-1}= b^{m-\delta_2} + b^{m-1}-2 b^{y_2},
\end{equation}
or, if $(a,b)=(2,5)$ and $(i,j,k,\kappa)=(3,4,2,1)$,
\begin{equation} \label{mum2}
2^{n-1}=5^{m-\delta_2}+2 \cdot 5^{m-1} -3 \cdot 5^{y_3}.
\end{equation}
Equation (\ref{mum}) has no solutions modulo $12$ if 
$$
(a,b) \in \{ (2,7), (3,5) \}.
$$ 
 If $(a,b)=(4,5)$, (\ref{gaew3}) implies that
 $$
5^{m-1}-5^{y_2} =  4^{x_4} + 5^{m}-4^n-5^{m-1},
 $$
 a contradiction modulo $3$. If $(a,b)=(2,5)$ and $(i,j,k,\kappa)=(3,4,2,1)$,  then we have equation (\ref{mum}), whence 
 $y_2=0$ or $\delta_2=m$. In the second case,
$$
\frac{1}{4} (N+4D) < 2^{n-1} < 5^{m-1} < \frac{1}{5} (N+4D),
$$
a contradiction. In the first, if $\delta_2=1$, we have
$$
2^{n-2} = 5^{m-1}-1,
$$
again a contradiction. If $\delta_2 \geq 2$, 
$$
\frac{1}{4} (N+4D) < 2^{n-1} \leq 5^{m-2}+5^{m-1}-2 < \frac{6}{25} (N+4D),
$$
once again a contradiction. If instead, we have equation (\ref{mum2}), so that $(i,j,k,\kappa)=(3,4,2,1)$, then $(\delta_1,\delta_2) = (1,1)$, and (\ref{gaew2}) becomes
$$
11 \cdot 5^{m-1} =  2^{n+1} - 3 \cdot 2^{x_2}.
$$
We thus have $x_2=0$ and hence a contradiction modulo $8$.
  
  We may thus suppose that $(a,b)=(2,5),$ that $(\delta_1,\delta_2)=(2,1)$ and that $(i,j,k,\kappa)=(2,3,4,1)$, whence (\ref{gaew2}) yields
  $$
  5^{m} = 11 \cdot 2^{n-2} - 2^{x_4+1},
  $$
  an immediate contradiction. 
  
  \subsection{Two $b$-large terms}

In case (\ref{hope2}), it follows that 
\begin{equation} \label{foosh!2}
D < \frac{1}{(j-i)a^{n-x_j}} (N+kD)
\end{equation}
whence
\begin{equation} \label{ook!2}
N+\kappa D > \left( \frac{(j-i)a^{n-x_j}+\kappa-k}{(j-i)a^{n-x_j}+4-k}. \right) (N+4D).
\end{equation}
Thus
$$
\left( \frac{1}{a^{\delta_1}} +\frac{1}{b^{\delta_2}} \right) (N+4D) > 
a^{n-\delta_1}+b^{m-\delta_2} = N+\kappa D >\left( \frac{(j-i)a^{n-x_j}+\kappa-k}{(j-i)a^{n-x_j}+4-k} \right) (N+4D)
$$
and so
\begin{equation} \label{spooky2}
 \frac{1}{a^{\delta_1}} +\frac{1}{b^{\delta_2}} >  \frac{(j-i)a^{n-x_j}+\kappa-k}{(j-i)a^{n-x_j}+4-k}.
 \end{equation}
 The identities
 $$
 (\kappa-j) (N+iD)+(j-i)(N+\kappa D) +(i-\kappa)(N+jD)=0
 $$
 and
 $$
 (\kappa-j) (N+kD)+(j-k)(N+\kappa D) +(k-\kappa)(N+jD)=0
 $$
 now lead  to the equations
  \begin{equation} \label{gaewA}
\left(  (j-i) b^{\delta_2}   + i-j \right) b^{m-\delta_2} = (\kappa -j) a^{x_i} + (j-i) a^{n-\delta_1} + (i-\kappa ) a^{x_j}
 \end{equation} 
 and
   \begin{equation} \label{gaewB}
\left(  (\kappa -j) a^{\delta_1}   + j-k \right) a^{n-\delta_1} +(k-\kappa) a^{x_j} = (j-\kappa ) b^{y_k} + (k-j) b^{m-\delta_2} + (\kappa-k ) b^m,
 \end{equation} 
 respectively.

 From (\ref{spooky2}), we find that $a \leq 4$.  If $a=4$, then, from (\ref{awesome}) and (\ref{spooky2}),  $x_j=n-1$, $\delta_1=1$ and $(i,j,k,\kappa)=(2,3,4,1)$.  It follows, if $a=4$, that
 $$
 3D=3 \cdot (4^{n-1}-4^{x_2}) = (N+4D)-(N+D) = 3 \cdot 4^{n-1} + b^{y_4}-b^{m-\delta_2},
 $$
 whence
 $$
 b^{m-\delta_2} - b^{y_4} =3 \cdot 4^{x_2}.
 $$
Lemma \ref{tech-lem} thus implies that either
$$
b=3 \cdot 4^{x_2}+1, \; m-\delta_2=1 \; \mbox{ and } \; y_4=0,
$$
or that 
$$
b=7, \;  m-\delta_2=2, \; y_4=0 \; \mbox{ and } \; x_2=2.
$$
 In the first case, we have
 $$
 N = (N+D)-D = 4^{x_2+1}+1,
 $$
 and so, from $N+(N+4D)=2(N+2D)$,
 $$
 b^m=2^{2n-1}+2^{2x_2}+1.
 $$
 From Theorem 1 of \cite{BeBuMi}, $b > 4$ and $m \geq 2$, it follows that 
 $$
 (b,m) \in \{ (7,2), (23,2) \},
 $$
 contradicting  $b=3 \cdot 4^{x_2}+1$. In the second case, we have
 $$
 D= \frac{1}{2} \left( 7^m-7^2 \right) = 4^{n-1}-4^{x_2}
 $$
 and so
 $$
 7^m-7^2=2^{2n-1}-2^{2x_2+1}.
 $$
 Proposition \ref{DT} leads to the desired contradiction.
 
 If $a=3$, we again have, from (\ref{awesome}) and (\ref{spooky2}),  that $x_j=n-1$, $\delta_1=1$, and $(i,j,k,\kappa)=(2,3,4,1)$ or $(3,4,2,1)$. If  $(i,j,k,\kappa)=(2,3,4,1)$, then, from (\ref{gaewB}),
 \begin{equation} \label{maniac}
 4 \cdot 3^{n-1} = 3 \cdot b^m - b^{y_4} - b^{m-\delta_2}
 \end{equation}
 and so $b$ is necessarily even. Equation (\ref{gaewA}) implies that
 $$
 (b^{\delta_2}-1) \, b^{m-\delta_2} = 2 \cdot (3^{n-1}-3^{x_2})
 $$
 and so, via parity considerations,  $m \neq \delta_2$. Thus $b \mid  4 \cdot 3^{n-1}$ and we can write $b=2^{\kappa_1} 3^{\kappa_2}$, where $\kappa_1 \in \{ 1, 2 \}$.  Equation (\ref{maniac}) thus reduces to (\ref{DT1}) with $(p,q)=(2,3)$ and hence we may appeal to Proposition \ref{DT} to conclude, after a little work, that there are no progressions corresponding to this case. If $a=3$ and $(i,j,k,\kappa)=(3,4,2,1)$, then we have
 $$
 D = 3^{n-1}-3^{x_3} = (N+2D)-(N+D) = 2 \cdot 3^{n-1} + b^{y_2}-b^{m-\delta_2},
 $$
 so that
 $$
\frac{1}{9} (N+4D) <  3^{n-1} +3^{x_3}+ b^{y_2}=b^{m-\delta_2} < 
b^{-\delta_2} (N+4D),
 $$
 whence
 $$
 4 \leq b \leq 8 \; \mbox{ and } \; \delta_2=1.
 $$
 From (\ref{gaewA}), we have
 $$
 \left(  b  -1 \right) b^{m-1} =   3^{n} -3^{x_3+1}
 $$
 and hence, applying Proposition \ref{DW}, after a short check, we conclude as desired.
 
 We therefore may suppose that $a=2$. Our goal is to show, in this case, that we have only $5$-term arithmetic progressions corresponding to
 $$
(a,b,N,D) \in \{ (2,3,17,24), (2,3,41,24), (2,9,17,24), (2,9,41,24) \}.
$$
 Let us assume first that $\delta_1=1$ and that $x_j=n-1$, whence (\ref{awesome}) and (\ref{spooky2}) imply that 
 $$
 (i,j,k,\kappa) \in 
 \{ (2,3,4,1), (2,4,3,1), (3,4,2,1), (1,3,4,2) \}.
 $$
 We will treat each of these four cases in turn. Suppose first that $(i,j,k,\kappa)=(2,3,4,1)$. Then 
 $$
 \left(  b^{\delta_2}   -1 \right) b^{m-\delta_2} = 2^n-2^{x_2+1} 
 $$
 and
 $$
2^n = 3 \cdot b^m - 2 b^{y_4} - b^{m-\delta_2},
 $$
 from (\ref{gaewA}) and (\ref{gaewB}), respectively. The second of these implies that
 $$
 2^n \geq 3 (b^m-b^{m-1}) > b^m
 $$
 and hence we can apply Theorem 1.6 of \cite{Be-Pill} to the first equation to conclude that $b=3$ and that
 $$
 (n,m,x_2,\delta_2) \in \{ (5,3,2,2), (8,5,3,4) \}.
 $$
 Neither coincides with a $5$-term progression (in each case $N+4D \not\in S_{2,3}$).  If $(i,j,k,\kappa)=(2,4,3,1)$ or $(3,4,2,1)$, from the identity $(N+D)+(N+4D)=(N+2D)+(N+3D)$, we find that
 $$
 b^{m-\delta_2} - b^{y_k}  = 2^{x_i},
 $$
 and so Lemma \ref{tech-lem} implies that
 we have one of
 \begin{equation} \label{ann1}
 b=2^{x_i}+1, \; m-\delta_2=1, \; y_k=0,
\end{equation}
\begin{equation} \label{ann3}
b=3, \; m-\delta_2=2, \; y_k=0, \;  x_i=3,
\end{equation}
or
\begin{equation} \label{ann4}
b=9, \; m-\delta_2=1, \; y_k=0, \;  x_i=3.
\end{equation}
In the first case, if $(i,j,k,\kappa)=(2,4,3,1)$, we have
$$
D = 2^{n-2}-2^{x_2-1} = \frac{1}{3} 2^{x_2} ( 2^{x_2}+1)^{m-1},
$$
a contradiction modulo $2^{x_2}$. If, on the other hand, we have (\ref{ann1}) and  $(i,j,k,\kappa)=(3,4,2,1)$,
$$
D = 2^{n-1}-2^{x_3} = \frac{1}{3} 2^{x_3} ( 2^{x_3}+1)^{m-1},
$$
so that
$$
2^{n-1-x_3}-1 = \frac{1}{3} ( 2^{x_3}+1)^{m-1}.
$$
Modulo $3$, $x_3$ is odd and $m \geq 2$. Thus, modulo $8$, we find that either $x_3=1$, $m=2$, $n=3$, or $x_3=1$, $m=3$, $n=4$, or $x_3=3$, $m=2$, $n=4$. Only the first of these corresponds to a $5$-term arithmetic progression in $S_{a,b}$, namely $(a,b,N,D)=(2,3,5,2)$, which we have previously encountered. If $(i,j,k,\kappa)=(2,4,3,1)$, (\ref{gaewA}) becomes
$$
\left(  b^{\delta_2}   -1 \right) b^{m-\delta_2} = 3 \cdot  (2^{n-2}-2^{x_2-1}),
$$
while $(i,j,k,\kappa)=(3,4,2,1)$ yields
$$
\left(   b^{\delta_2}   -1 \right) b^{m-\delta_2} = 3 \cdot (2^{n-1} - 2^{x_3}).
$$
In cases (\ref{ann3}) and  (\ref{ann4}), we may thus apply Proposition \ref{DT} to these equations, concluding after a short computation  with the $5$-term arithmetic progressions corresponding to
$$
(a,b,N,D) = (2,3,17,24) \; \mbox{ and } \; (2,9,17,24).
$$
Finally,   
 if $(i,j,k,\kappa)=(1,3,4,2)$, then, from $(N+D)+(N+4D)=(N+2D)+(N+3D)$, we have
  $$
  b^{m-\delta_2}-b^{y_4}=2^{x_1},
  $$
  while (\ref{gaewB}) gives
  $$
 2^{n-1}  = 2 b^m -b^{y_4} -  b^{m-\delta_2}.
  $$
  Lemma \ref{tech-lem} thus implies that $y_4=0$ and so combining the preceding two equations,
  $$
  b^m =2^{n-2} +2^{x_1-1}+1.
  $$
  Since necessarily $m \geq 2$ (else we have three equal $y_i$), we may appeal to Theorem 1 of \cite{BeBuMi} to conclude that
  $$
  b^m \in \{ 7^2, 23^2, 3^4 \},
  $$
  so that 
  $$
  (b,m,n,x_1) \in \{ (7,2,7,5), (23,2,11,5), (3,4,8,5), (9,2,8,5) \};
  $$
  none of these correspond to a $5$-term arithmetic progression in $S_{a,b}$.

Continuing under the assumption that $\delta_1=1$, let us next suppose that $n-x_j=2$. From (\ref{awesome}) and (\ref{spooky2}), we have that either
 $$
 (i,j,k,\kappa) \in 
 \{ (2,3,4,1), (3,4,2,1) \},
 $$
 or that
 \begin{equation} \label{speccy}
 (b,\delta_2,i,j,k,\kappa) \in \{  (3,1,2,4,3,1), (4,1,2,4,3,1), (5,1,2,4,3,1), (3,1,1,3,4,1) 
 \}.
 \end{equation}
  Assume first that $(i,j,k,\kappa) = 
 (2,3,4,1)$. Then we have
 $$
 N = (N+4D)-4D = 2^{x_2+2}+b^{y_4}.
 $$
 It follows from the identity  $N+(N+3D)=(N+D)+(N+2D)$  that 
 \begin{equation}\label{from}
  3 \cdot 2^{x_2}+b^{y_4}=2^{n-2}+b^{m-\delta_2},
 \end{equation}
 while (\ref{gaewA}) gives
  \begin{equation}\label{from2}
 \left(  b^{\delta_2}   -1 \right) b^{m-\delta_2} = 3 \cdot 2^{n-2}- 2^{x_2+1}. 
  \end{equation}
 If $x_2=n-3$, (\ref{from}) becomes
 $$
 2^{n-3}=b^{m-\delta_2}-b^{y_4},
 $$
 while (\ref{from2}) yields
 $$
 2^{n-1} = b^m- b^{m-\delta_2}.
  $$
Appealing to Lemma \ref{tech-lem}, the second equation implies that $m-\delta_2=0$, contradicting the first.  If $x_2=n-4$, (\ref{from}) becomes
 $$
b^{y_4}-b^{m-\delta_2}=2^{n-4},
 $$
 and hence Lemma \ref{tech-lem} implies that $m-\delta_2=0$,
 and that 
 $$
 (b,y_4,n) \in \{ (2^{n-4}+1,1,n), (3,1,5), (3,2,7), (9,1,7) \}.
 $$
 Since (\ref{from2}) yields
 $$
  \left(  b^{\delta_2}   -1 \right) b^{m-\delta_2} = 5 \cdot 2^{n-3},
  $$
we have $b^m-1=5 \cdot 2^{n-3}$, whence $ (b,y_4,n) \in \{ (3,2,7), (9,1,7) \}$ correspond to the progressions 
 $(a,b,N,D) \in \{  (2,3,41,24), (2,9,41,24) \}$ and $(b,y_4,n)= (3,1,5)$ does not lead to a progression. If $(b,y_4)=(2^{n-4}+1,1)$, 
 we have
 \begin{equation} \label{kook}
 \left( 2^{n-4}+1 \right)^m-1 = 5 \cdot 2^{n-3}.
 \end{equation}
The values $n \in \{ 4, 5, 6 \}$ lead to no progressions, while $n=7$ corresponds to $(a,b,N,D)=(2,9,41,24)$. If $n \geq 8$, equation (\ref{kook}) contradicts $m \geq 2$.
 
 Next, suppose that we have $a=2$, $\delta_1=1$, $n-x_j=2$ and 
$(i,j,k,\kappa) =  (3,4,2,1)$. From $D=(N+2D)-(N+D)$, we have
\begin{equation} \label{fossil}
2^{n-2}+2^{x_3}=b^{m-\delta_2}-b^{y_2}
\end{equation}
and so
$$
\frac{1}{8} (N+4D) < 2^{n-2} < b^{m-\delta_2}-b^{y_2} < b^{-\delta_2} (N+4D).
$$
Thus $\delta_2=1$ and $b \leq 7$. From $3D=(N+4D)-(N+D)$, we further have
$$
 (b-1) b^{m-1} = 2^n - 3 \cdot 2^{x_3},
 $$
 an immediate contradiction modulo $3$ for $b \in \{ 3, 4, 6, 7 \}$. If $b=5$, applying Proposition \ref{DT} to equation (\ref{fossil}), we find no new arithmetic progressions. 
 
 If $a=2$, $\delta_1=1$, $n-x_j=2$ and we have (\ref{speccy}), then, via (\ref{gaewB}), we find that either
 $$
  0=3^{y_4}- 5 \cdot 3^{m-1},
  $$
  or that 
  $$
  -2^n=3 \cdot b^{y_3}-(2b+1) \cdot b^m.
  $$
  Only the case $(b,i,j,k,\kappa) =(5, 2, 4, 3, 1)$ is not an immediate contradiction, where we have
  $$
  -2^n=3 \cdot 5^{y_3}-11 \cdot 5^m.
  $$
  From Proposition \ref{DW}, we find that $m=0$, a contradiction.
 
 If $a=2$, $\delta_1=1$ and  $n-x_j=3$, then (\ref{awesome}) and (\ref{spooky2}) imply that  $\delta_2=1$ and that either
 $$
 3 \leq b \leq 7 \; \mbox{ and } \; (i,j,k,\kappa)=(2,3,4,1),
 $$
 or 
  \begin{equation} \label{speccy2}
 (b,i,j,k,\kappa) \in \{  (3,2,4,3,1), (3,3,4,2,1), (4,3,4,2,1) \}.
  \end{equation}
 If $(i,j,k,\kappa)=(2,3,4,1)$, (\ref{gaewB}) becomes
  \begin{equation} \label{muck1}
  -5 \cdot 2^{n-3} = 2 \cdot b^{y_4} +(1-3b) \cdot b^{m-1},
  \end{equation}
  while (\ref{gaewA}) yields
  \begin{equation} \label{muck2}
  \left(  b -1 \right) \cdot  b^{m-1} =- 2^{x_2+1} + 5 \cdot 2^{n-3},
   \end{equation}
  contradictions modulo $5$ for $5 \leq b \leq 7$. If $b=4$, then $n=2m$ or $n=2m+1$. In the first case, (\ref{muck1}) implies that $2^{2m-2} \, \| \, 2 \cdot b^{y_4}$, contradicting $b=4$. In the second, $3 \cdot 2^{2m-1} = 2  \cdot 4^{y_4}$, again a contradiction. If $b=3$, (\ref{muck1}) implies that $n=4$ and $y_4=0$, a contradiction modulo $8$.  In case (\ref{speccy2}), 
 we find from (\ref{gaewB}) that
 $$
2^{n-2} = 7 \cdot 3^{m-3} - 3^{y_k-1},
 $$
 $$
15 \cdot 2^{n-3} =5 \cdot  3^{m-1} -3^{y_k+1} 
 $$
 and
 $$
5 \cdot 2^{n-3} =2^{2m-1} - 2^{2y_k},
 $$
 respectively. The first two equations lead to contradictions modulo $24$ and $5$, while the third implies that $5 \cdot 2^{n-3} <2^{2m-1}$, contradicting $n=2m$ or $n=2m+1$.  
  
  If we have $\delta_1=1$ and $n-x_j \geq 4$, (\ref{awesome}) and (\ref{spooky2}) imply that $n-x_j = 4$, $\delta_2=1$ and that
  $$
  (b,i,j,k,\kappa)=(3,2,3,4,1),
  $$
  whence, from (\ref{gaewB}), 
  $$
5 \cdot 2^{n-5} = 4 \cdot 3^{m-1} -3^{y_k},
  $$
  and so $n=5$, an immediate contradiction.
  
  We may therefore assume, for the remainder of this section, that $\delta_1 \geq 2$. From (\ref{awesome}) and (\ref{spooky2}), $a=2$, and either $(i,j,k,\kappa)=(2,3,4,1)$, or  we have $\delta_2=1$, and one of
  \begin{equation} \label{speccy3}
  \delta_1=2, \; 3 \leq b \leq 6, \; n-x_j=1 \;  \mbox{ and } \; (i,j,k,\kappa)=(2,4,3,1), 
  \end{equation}
  \begin{equation} \label{speccy4}
  \delta_1=2, \; 3 \leq b \leq 5, \; n-x_j=1 \; \mbox{ and } \; (i,j,k,\kappa)=(3,4,2,1), 
  \end{equation}
  \begin{equation} \label{speccy5}
  \delta_1=2, \; b=3, \; n-x_j=2, \;  (i,j,k,\kappa) = (3,4,2,1),
   \end{equation}
  \begin{equation} \label{speccy6}
  \delta_1=3, \; b=3, \; n-x_j=1, \;  (i,j,k,\kappa) \in \{ (2, 4, 3, 1), (3,4,2,1) \},
   \end{equation}
  or
  \begin{equation} \label{speccy7}
  \delta_1 \geq 4, \; b=3, \; n-x_j=1, \;  (i,j,k,\kappa) =(3,4,2,1).
   \end{equation}
   
 If $(i,j,k,\kappa)=(2,3,4,1)$ and $\delta_1=2$, then $n-x_3 \in \{ 1, 2 \}$ and  (\ref{gaewB}) implies
$$
3 (2^{n-1}-2^{x_3})  =3 b^m - 2 b^{y_4} -  b^{m-\delta_2}. 
$$
Thus either
$$
0 =3 b^m - 2 b^{y_4} -  b^{m-\delta_2},
$$
in case $n-x_3=1$, or
\begin{equation} \label{ooo!}
3  \cdot 2^{n-2}  =3 b^m - 2 b^{y_4} -  b^{m-\delta_2},
\end{equation}
if $n-x_3=2$.
The first of these contradicts $y_4<m$. In case (\ref{ooo!}), either $y_4=0$ or $\delta_2=m$, or we have that $b=2^\alpha 3^\beta$, where $\beta \in \{ 0, 1 \}$. 
If $y_4=0$, 
$$
D =3 \cdot 2^{n-2}+1-b^m=2^{n-2} -2^{x_2},
$$
so
$$
b^m = 2^{n-1}+2^{x_2}+1.
$$
By Theorem 1 of \cite{BeBuMi}, this contradicts $m \geq 2$, unless 
$$
(b,m,n,x_2) \in \{ (7,2,6,4), (23,2,10,4), (3,4,7,4), (9,2,7,4) \}.
$$
None of these correspond to $5$-term progressions in $S_{2,b}$.
If $\delta_2=m$ and $n-x_3=2$, then the identity $4(N+D)=(N+4D)+3N$ implies that
$$
b^{y_1} = b^{y_4} + 3 \cdot 2^{x_2}.
$$
But then
$$
D=(N+4D)-(N+3D)=3 \cdot 2^{n-2}+b^{y_4}-b^m =b^{y_1}-b^m < 0.
$$
We may thus suppose that $b=2^\alpha 3^\beta$, where $\beta \in \{ 0, 1 \}$ and $y_4(m-\delta_2) \neq 0$. Since (\ref{gaewA}) implies that
\begin{equation} \label{skoot}
\left(  b^{\delta_2}   -1 \right) b^{m-\delta_2} = 2^{n-1}-2^{x_2+1},
\end{equation}
we have $(m-\delta_2)\alpha = x_2+1$ and so
$$
\left( 2^{\delta_2 \alpha} 3^{\delta_2\beta}-1 \right) 3^{(m-\delta_2)\beta} = 2^{n-x_2-2}-1.
$$
If $\beta=0$, we have $m\alpha = n-1$ and so (\ref{ooo!}) implies that
$$
3  \cdot 2^{n-2}  = 2^{\alpha y_4+1} +  2^{(m-\delta_2)\alpha},
$$
so that, from $(m-\delta_2)\alpha < m\alpha = n-1$,
$$
(\alpha y_4+1,(m-\delta_2)\alpha) = (n-1,n-2).
$$
But then $\delta_2\alpha =1$, contradicting $b > 2$. We may thus suppose that $\beta=1$, so that, from (\ref{ooo!}), either $y_4=1$ or $m-\delta_2=1$. In the latter case, from (\ref{skoot}), 
$$
2^{\delta_2 \alpha-1} 3^{\delta_2+1} = 2^{n-x_2-3}+1,
$$
so that $\delta_2= \alpha=1$, $b=6$, $n-x_2=6$ and (\ref{ooo!}) becomes
$$
2^{n-3}  =17 \cdot 6^{m-2} - 2 \cdot 6^{y_4-1},
$$
whence $y_4=1$, a contradiction modulo $4$. Suppose then that $y_4=1$ and $m-\delta_2 \geq 2$, so that, from (\ref{ooo!}),
\begin{equation} \label{earl}
2^{n-2}  =3^{m} 2^{\alpha m} - 2^{1+\alpha}   -  3^{m-\delta_2-1} 2^{\alpha(m-\delta_2)}.
\end{equation}
If $\alpha \geq 2$, it follows that $n=\alpha+3$, whence
$$
N+4D > b^m \geq b^2 =9 \cdot 2^{2\alpha} > 2^{n+1} > N+4D.
$$
Thus $\alpha \in \{ 0,1 \}$. If $\alpha=0$, (\ref{earl}) yields
$$
2^{n-2} +2 =3^{m} -3^{m-\delta_2-1},
$$
and so Proposition \ref{DT} and a short computation reveals so solutions with $\delta_2 \geq 1$.
If $\alpha=1$, we have
$$
2^{n-2}  =6^{m} - 4   -  3^{m-\delta_2-1} 2^{m-\delta_2}
$$
and so $m-\delta_2=2$ and 
$$
2^{n-4}  =2^{m-2} 3^{m} - 1   -  3^{m-\delta_2-1},
$$
a contradiction modulo $8$, for $m \geq 5$. A short check finds no new $5$-term progressions with $(a,b)=(2,6)$ and $m \leq 4$.

Next, suppose that $(i,j,k,\kappa)=(2,3,4,1)$ and that $\delta_1 \geq 3$. Then, from (\ref{awesome}) and (\ref{spooky2}), $\delta_1=3$, $\delta_2=1$ and $n-x_3 \in \{ 1, 2 \}$.
Equation (\ref{gaewA}) and (\ref{gaewB}) now become
\begin{equation} \label{tool1}
\left(   b -1 \right) b^{m-1} = 5 \cdot 2^{n-3} -2^{x_2+1} 
\end{equation} 
 and
\begin{equation} \label{tool2}
5 \cdot  2^{n-3}  = 3 \cdot b^m - 2 \cdot b^{y_4} -  b^{m-1},
\end{equation}
or
\begin{equation} \label{tool3}
\left( b -1 \right) b^{m-1} = 3 \cdot 2^{n-3} -2^{x_2+1} 
\end{equation}
 and
\begin{equation} \label{tool4}
11 \cdot  2^{n-3} = 3 \cdot b^m - 2 \cdot b^{y_4} -  b^{m-1},
\end{equation}
if $n-x_3=1$ or $2$, respectively. Combining (\ref{tool1}) and (\ref{tool2}), we have
$$
b^m - b^{y_4} =2^{x_2},
$$
whence, via Lemma \ref{tech-lem}, $b=3$ and $m=2$, an easy contradiction.
Conversely, if we have (\ref{tool3}) and (\ref{tool4}), 
$$
b^m - 4 \cdot b^{m-1}  = - 3 \cdot b^{y_4}  - 11 \cdot 2^{x_2},
$$
and so $b=3$ and 
$$
3^{m-1}  = 3^{y_4+1}  + 11 \cdot 2^{x_2},
$$
contradicting $m \geq 2$.

We may thus suppose that $\delta_2=1$ and that one of (\ref{speccy3})--(\ref{speccy7}) holds. If $(i,j,k,\kappa)=(2,4,3,1)$, then we are in cases (\ref{speccy3}) or (\ref{speccy6}), and (\ref{gaewA}) and (\ref{gaewB}) yield
$$
\left(   b  -1 \right) b^{m-1} = -3 \cdot 2^{x_2-1} + 2^{n-\delta_1} + 2^{n-2}
$$
 and
$$
\left(  -3 \cdot 2^{\delta_1}   +1 \right) 2^{n-\delta_1} + 2^{n} = 3 \cdot b^{y_3} - b^{m-1} -2 \cdot b^m.
$$
If $\delta_1=3$, then $b=3$ and these equations become
$$
3^{m-2} = 2^{n-4}-2^{x_2-2} 
$$
 and
$$
5 \cdot 2^{n-3}  = 7 \cdot 3^{m-2} - 3^{y_3};
$$
the latter equation is a contradiction modulo $8$ since we may suppose that $n \geq 6$. If $\delta_1=2$, we have
$$
\left(   b  -1 \right) b^{m-1} = 2^{n-1} -3 \cdot 2^{x_2-1} 
$$
 and
$$
7 \cdot 2^{n-2}  = 2 \cdot b^m-3 \cdot b^{y_3} + b^{m-1},
$$
where $3 \leq b \leq 6$. The first equation is a contradiction modulo $3$ unless $b=5$, in which case the second equation becomes
$$
7 \cdot 2^{n-2}  = 11 \cdot 5^{m-2}-3 \cdot 5^{y_3}.
$$
Proposition \ref{DW} and a short computation show that there are. no solutions.

If $(i,j,k,\kappa)=(3,4,2,1)$ and we have one of (\ref{speccy4})--(\ref{speccy7}), then, from (\ref{gaewA}) and (\ref{gaewB}), 
\begin{equation} \label{cook1}
\left(  b   -1 \right) b^{m-1} = -3 \cdot 2^{x_3} + 2^{n-\delta_1} +2^{x_4+1}
\end{equation}
 and
\begin{equation} \label{cook2}
\left(  -3 \cdot 2^{\delta_1}   +2 \right) 2^{n-\delta_1} + 2^{x_4} = 3 \cdot b^{y_2} -2 \cdot b^{m-1} - b^m.
\end{equation}
In case  (\ref{speccy4}), we have 
$$
\left(  b -1 \right) b^{m-1} = 5 \cdot 2^{n-2}-3 \cdot 2^{x_3} 
$$
 and
$$
2^{n+1} = b^m +2 \cdot b^{m-1} - 3 \cdot b^{y_2},
$$
with $3 \leq b \leq 5$; the first equation leads to an immediate contradiction modulo $3$ or $5$. In case (\ref{speccy5}), $b=3$ and $\delta_1=n-x_4=2$, whence (\ref{cook1}) yields
$$
3^{m-2} = 2^{n-3} -2^{x_3-1}.
$$
 It follows that $x_3=1$ and that $(m,n)=(4,6)$ or $(3,4)$; each contradicts $N \geq 2$.
  If we have (\ref{speccy6}),  (\ref{cook1}) implies that
  $$
 3^{m-2} = 3 \cdot 2^{n-4} -2^{x_3-1},
  $$
  contradicting $n \geq 5$. From (\ref{speccy7}) and (\ref{cook1}),
  $$
3^{m-1} = -3 \cdot 2^{x_3-1} + 2^{n-\delta_1-1} +2^{n-1};
  $$
 Appeal to Proposition \ref{DT} and a short computation contradicts $\delta_1 \geq 4$.

%---------------------------------------------------------------
\section{Two large terms} \label{Sec6}
%---------------------------------------------------------------

Next, suppose we have precisely two large terms. We require
$$
\frac{1}{2}(N+4D) < N+2D \leq a^{n-1}+b^{m-1} < \left( \frac{1}{a} +
\frac{1}{b} \right) (N+4D),
$$
so that $a=2$ or 
$$
(a,b) \in \{ (3,4), (3,5) \}.
$$
If $N+4D$ is both $a$-large and $b$-large, then
$$
\frac{1}{2}(N+4D) < N+2D \leq a^{n-1}+b^{m-1} <  \frac{1}{a} (N+4D),
$$
contradicting $a \geq 2$. We next show that
$N+3D$ is large.
If not then 
$$
\frac{3}{4}(N+4D) < N+3D \leq a^{n-1}+b^{m-1} < \left( \frac{1}{a} +
\frac{1}{b} \right) (N+4D),
$$
and so $(a,b)=(2,3)$ and  $N+3D= 2^{n-1}+3^{m-1}$. If $N+2D$ is also not large, then necessarily $x_2 = n-1$, and so, from $2(N+3D)=(N+2D)+(N+4D)$, either
$$
2^n+2 \cdot 3^{m-1}=2^{n-1}+3^{y_2}+2^n+3^{y_4}
$$
or
$$
2^n+2 \cdot 3^{m-1}=2^{n-1}+3^{y_2}+2^{x_4}+3^{m}.
$$
The first of these implies that
$$
2 \cdot 3^{m-1}=2^{n-1}+3^{y_2}+3^{y_4},
$$
while the second gives
$$
2^{n-1}=3^{y_2}+2^{x_4}+3^{m-1};
$$
in either case, we may appeal to Proposition \ref{DT} and a short computation to confirm that there are no corresponding new arithmetic progressions.  We may thus suppose that $N+2D$ is large, whence
we have one of
$$
2^n+2 \cdot 3^{m-1}=2^{n}+3^{y_2}+2^n+3^{y_4},
$$
$$
2^n+2 \cdot 3^{m-1}=2^{x_2}+3^{m}+2^n+3^{y_4},
$$
$$
2^n+2 \cdot 3^{m-1}=2^{n}+3^{y_2}+2^{x_4}+3^{m},
$$
or
$$
2^n+2 \cdot 3^{m-1}=2^{x_2}+3^{m}+2^{x_4}+3^{m}.
$$
These yield
$$
2 \cdot 3^{m-1}=3^{y_2}+2^n+3^{y_4},
$$
$$
1=2^{x_2}+3^{m-1}+3^{y_4},
$$
$$
1=3^{y_2}+2^{x_4}+3^{m-1},
$$
and
$$
2^n=2^{x_2}+4 \cdot 3^{m-1}+2^{x_4},
$$
respectively, again all treatable via appeal to Proposition \ref{DT}. Once more, we find no new $5$-term arithmetic progressions in $S_{2,3}$.

We may thus suppose that both $N+3D$ and $N+4D$ are large, and that $N+D$ and $N+2D$ are not large. If $x_2 \leq n-2$, then, from $2(N+2D) > N+4D$, necessarily $x_2=n-2, y_2=m-1$, $a=2$ and $b=3$. But then, we have one of
$$
2^{n+1} + 2 \cdot 3^{y_3}=2^{n-2}+3^{m-1}+2^{n}+3^{y_4},
$$
$$
2^{n+1} + 2 \cdot 3^{y_3}=2^{n-2}+3^{m-1}+2^{x_4}+3^{m},
$$
$$
2^{x_3+1} + 2 \cdot 3^{m}=2^{n-2}+3^{m-1}+2^{n}+3^{y_4},
$$
or
$$
2^{x_3+1} + 2 \cdot 3^{m}=2^{n-2}+3^{m-1}+2^{x_4}+3^{m}.
$$
These reduce to
$$
3 \cdot 2^{n-2} + 2 \cdot 3^{y_3}=3^{m-1}+3^{y_4},
$$
$$
2^{n+1} + 2 \cdot 3^{y_3}=2^{n-2}+4 \cdot 3^{m-1}+2^{x_4},
$$
$$
2^{x_3+1} + 5 \cdot 3^{m-1}=5 \cdot 2^{n-2}+3^{y_4},
$$
and
$$
2^{x_3+1} + 2 \cdot 3^{m-1}=2^{n-2}+2^{x_4}.
$$
We apply Proposition \ref{DT} to the first and fourth of these equations and Proposition \ref{BajBen} to the second; we find no new $5$-term progressions in $S_{2,3}$.
For the third equation, which corresponds to the case where $N+3D$ is $3$-large and $N+4D$ is $2$-large, we note that the 
identity $2(N+2D)=(N+D)+(N+3D)$ yields
$$
2^{n-1}= 2^{x_1}+3^{y_1}+2^{x_3}+3^{m-1}.
$$
Once again, Proposition \ref{BajBen} and a short calculation find no corresponding $5$-term progressions in $S_{2,3}$.

We may thus suppose that $x_2=n-1$. If 
$$
(a,b) \in \{ (3,4), (3,5) \},
$$
necessarily $y_2=m-1$ and, from $2(N+3D)=(N+2D)+(N+4D)$, 
\begin{equation} \label{c1}
2 \cdot 3^{n-1}+2 \cdot b^{y_3} = b^{y_2} +b^{y_4},
\end{equation}
\begin{equation} \label{c2}
5 \cdot 3^{n-1}+2 \cdot b^{y_3} = b^{y_2} + 3^{x_4}+b^{m},
\end{equation}
\begin{equation} \label{c3}
2 \cdot 3^{x_3}+2 \cdot b^{m} = 4 \cdot 3^{n-1}+b^{y_2} +b^{y_4},
\end{equation}
or
\begin{equation} \label{c4}
2 \cdot 3^{x_3}+b^{m} = 3^{n-1}+b^{y_2} + 3^{x_4}.
\end{equation}
In cases (\ref{c1}), (\ref{c3}) and (\ref{c4}), if $b=4$, we conclude as desired from Proposition \ref{BajBen}.  In case (\ref{c2}) and $b=4$, the identity $2(N+2D)=(N+D)+(N+3D)$ leads to the equation
$$
2 \cdot 4^{y_2} = 3^{x_1}+4^{y_1}+ 3^{n-1} + 4^{y_3};
$$
once again, Proposition \ref{BajBen}
leads to the desired conclusion.

We next consider the case $b=5$ in equations (\ref{c1})--(\ref{c4}). We have $N+2D=3^{n-1}+5^{m-1}$ and 
either $x_1=n-1$, or $(x_1,x_2)=(n-2,m-1)$. In the latter case, $2(N+2D)=(N+D)+(N+3D)$ yields
$$
5 \cdot 3^{n-2}+5^{m-1}= 3^{x_3}+5^{y_3},
$$
and Proposition \ref{DT} shows that there are no corresponding $5$-term arithmetic progressions in $S_{3,5}$.
If $x_1=n-1$, we have
$$
3^{n-1}+2 \cdot 5^{m-1}=5^{y_1}+3^{x_3}+5^{y_3}.
$$
If $x_3=n$, this equation becomes
$$
2 \cdot 5^{m-1}=2 \cdot 3^{n-1}+ 5^{y_1}+5^{y_3},
$$
a contradiction modulo $4$, while $y_3=m$ yields
$$
3^{n-1}=5^{y_1}+ 3^{x_3}+3 \cdot 5^{m-1}.
$$
Again Proposition \ref{DT} leads to the desired conclusion.
It remains, then,  to handle the cases with $a=2$ and $x_2=n-1$. 

\subsection{The case $x_1=n-1$}

If  $x_1=x_2=n-1$, then necessarily $x_0 \leq n-2$. 

\subsubsection{ The subcase $x_0=n-2$}

If $x_0=n-2$, then the 
identity $N+(N+2D)=2(N+D)$ gives
\begin{equation} \label{poop}
b^{y_0}+b^{y_2}=2^{n-2}+2 \cdot b^{y_1}.
\end{equation}
If two of $y_0, y_1$ or $y_2$ are equal (note that $y_1 < y_2$), then this equation reduces to 
$$
b^{y_2}-b^{y_1}=2^{n-3},  
$$
and so, from Lemma \ref{tech-lem}, $y_1=0$, and either $b=2^{n-3}+1$, $y_2=1$, or we have 
$$
(b,y_2,n) \in \{ (3,1,4), (3,2,6), (9,1,6) \}.
$$
A short check eliminates these last three cases, while $b=2^{n-3}+1$ implies that
$$
N+4D > \left( 2^{n-3}+1 \right)^m \geq \left( 2^{n-3}+1 \right)^2 > 2^{n+1} > N+4D, 
$$
a contradiction, provided $n \geq 7$.
Thus $y_0 \neq y_2$ and $y_0 \neq y_1$,  and so 
$$
\frac{1}{8} (N+4D) < 2^{n-2} < \left( \frac{1}{b} + \frac{1}{b^2} \right) (N+4D),
$$
whence $b \leq 7$.  It thus follows from (\ref{poop}) that either $y_1y_0=0$, or that $b=4$. In $b \in \{ 3, 5, 6, 7 \}$, we can appeal to Proposition \ref{DT} to bound $n$ and hence $m$
; we find no new arithmetic progressions. If $b=4$, then, since $y_0 \neq y_2$ and $y_0 \neq y_1$, (\ref{poop}) and the fact that $y_2>y_1$ imply that $y_0y_1 \neq 0$, and that $n=2y_2+2$ and $2y_0=2y_1+1$, the latter an immediate contradiction. 

\subsubsection{ The subcases $x_0 \leq n-3$}

Suppose then that  $x_1=x_2=n-1$ and $x_0 \leq n-3$. Notice first that if also $x_4 \leq n-3$, we have
$$
N = 2^{x_0}+b^{y_0} < \left( \frac{1}{8} + \frac{1}{b} \right) (N+4D)
$$
and
$$
4D =2^{x_4+2}-2^{x_3+2} < 2^{n-1} < \frac{1}{2}  (N+4D),
$$
so
$$
N+4D < \left( \frac{5}{8} + \frac{1}{b} \right) (N+4D),
$$
contradicting $b \geq 3$. We may therefore assume that $x_4=n$ or $x_4=n-2$.  If  $x_3=x_4=n$, it follows from the identity $(N+D)+(N+4D)=(N+2D)+(N+3D)$ that
$$
b^{y_1}+b^{y_4} = b^{y_2}+b^{y_3},
$$
where $y_1<y_2$ and $y_3<y_4$. Thus $y_1=y_3$ and $y_2=y_4$, contradicting the fact that $y_3< y_4<m$. If, on the other hand, $y_4=m$, so that $x_4=n-2$, 
 from $3(N+2D)=(N+4D)+2(N+D)$, we have
\begin{equation} \label{quest}
2^{n-2} + 3 \cdot b^{y_2} = b^m +  2 \cdot b^{y_1},
\end{equation}
and so
$$
N+4D = 2^{n-2}+b^m = 2^{n-1} + 3 \cdot b^{y_2} -2 \cdot b^{y_1}.
$$
If $y_2 \leq m-2$, 
$$
N+4D  = 2^{n-1} + 3 \cdot b^{y_2} -2 \cdot b^{y_1} <
\left( \frac{1}{2} + \frac{3}{b^2} \right) (N+4D),
$$
contradicting $b \geq 3$. We thus have $y_2=m-1$, so that 
$$
N+4D =  2^{n-1} + 3 \cdot b^{m-1} -2 \cdot b^{y_1} <
\left( \frac{1}{2} + \frac{3}{b} \right) (N+4D),
$$
whence $b \leq 5$. The case $b=4$ contradicts (\ref{quest}) modulo $3$. If $b \in \{3, 5 \}$, (\ref{quest}) becomes
$$
2^{n-2} = 2 \cdot 3^{y_1}
$$
and
$$
2^{n-3}  =5^{m-1} + 5^{y_1},
$$
respectively, both contradictions modulo $4$. 

We are left, then, to treat the case with $x_4=n$ and $y_3=m$. 
The identity $3(N+2D)=(N+4D)+2(N+D)$ now implies 
\begin{equation} \label{quest2}
3 \cdot b^{y_2} = 2^{n-1}+ b^{y_4} +  2 \cdot b^{y_1},
\end{equation}
and hence
$$
N+4D=2^n+b^{y_4} < 6 \cdot b^{y_2}.
$$
We thus have, again, that $y_2=m-1$ and that $b \leq 5$. The assumption that $b=4$ contradicts (\ref{quest2}) modulo $3$. If $b=5$ and $y_0 \leq m-2$, then
$$
N+4D \leq 2^{n-3}+5^{m-2}+4(5^{m-1}-5^{y_1})=
2^{n-3}+21 \cdot 5^{m-2}-4 \cdot 5^{y_1} < \left( \frac{1}{8} + \frac{21}{25} \right) (N+4D).
$$
Thus $y_0=m-1$ and hence, from $N+(N+4D)=2(N+2D)$,
$$
5^{m-1}=5^{y_4}+2^{x_0},
$$
whence $m=2$, $y_4=0$, $x_0=2$, contradicting $N > 1$.
 If $b=3$, $N+(N+4D)=2(N+2D)$ implies that
$$
2 \cdot 3^{m-1}=2^{x_0}+3^{y_0}+3^{y_4},
$$
whereby $y_0y_4=0$. If $y_0=m-1$ or $y_4=m-1$, we  have
$$
3^{m-1}=2^{x_0}+1,
$$
a contradiction. Thus $3^{y_0}+3^{y_4} \leq 3^{m-2}+1$ and so
$$
2 \cdot 3^{m-1} \leq 2^{n-3}+3^{m-2}+1,
$$
i.e.
$$
5 \cdot 3^{m-2} \leq 2^{n-3}+1. 
$$
But then
$$
N+4D < 3^{m+1} \leq \frac{27}{5} \left(2^{n-3}+1 \right) 
< \frac{27}{40} (N+4D) +6,
$$
again a contradiction.

\subsection{The case $x_1=n-2$}

Next, let us suppose that $x_1=n-2$, so that
\begin{equation} \label{fonz}
N=2(N+D)-(N+2D) = 2 \cdot b^{y_1} - b^{y_2} > 1,
\end{equation}
and hence $y_1 \geq y_2$. If also $x_4=n$, $3(N+2D)=(N+4D)+2(N+D)$ implies that 
$$
3 \cdot b^{y_2}=b^{y_4} + 2 \cdot b^{y_1},
$$
whereby, from $y_1 \geq y_2$, 
we reach a contradiction (since we do not have $y_1=y_2=y_4$). We may thus assume that $y_4=m$. Let us first suppose that also $y_3=m$. If $x_4 \leq n-3$, then $D < 2^{n-3}$ and so
$$
N+4D = 2^{n-2}+b^{y_1}+ 3(2^{x_4}-2^{y_3}) < 5 \cdot 2^{n-3}+b^{y_1} < \left( \frac{5}{8} + \frac{1}{b} \right) (N+4D),
$$
contradicting $b \geq 3$.
Thus $x_4 \in \{ n-2, n-1 \}$. If $x_4=n-2$, from $(N+D)+(N+4D)=(N+2D)+(N+3D)$, we have
$$
b^{y_1} = b^{y_2}+2^{x_3},
$$
and so, from Lemma \ref{tech-lem}, $y_2=0$ and either
$$
b=2^{x_3}+1, \; y_1=1
$$
or
$$
(b,y_1,x_3) \in \{ (3,1,1), (3,2,3), (9,1,3) \}.
$$
Since the identity $3(N+2D)=(N+4D)+2(N+D)$ implies
\begin{equation} \label{pain}
3 \cdot 2^{n-2}= b^m-2 \cdot b^{y_1}-3,
\end{equation}
we have a contradiction modulo $8$, unless 
$$
(b,y_1,x_3) \in \{ (3,1,1), (5,1,2) \}.
$$
With (\ref{pain}), we thus have
$$
2^{n-2}= 3^{m-1}-3
$$
or
$$
3 \cdot 2^{n-2}= 5^m-13.
$$
The first is an immediate contradiction, while the second, with Proposition \ref{DW}, implies that $(n,m)=(4,2)$, which fails to correspond to a new progression. We thus have $x_4=n-1$ and so, again from $(N+D)+(N+4D)=(N+2D)+(N+3D)$, 
\begin{equation} \label{colic}
2^{n-2}+b^{y_1}=2^{x_3} + b^{y_2}.
\end{equation}
From $y_1 \geq y_2$ and (\ref{colic}), either $y_1=y_2$ and $x_3=n-2$, or we have $x_3 =n-1$, whereby
$$
b^{y_1}-b^{y_2}=2^{n-2}.
$$
In this latter case, from Lemma \ref{tech-lem} and $n \geq 6$, $y_2=0$, $y_1=1$ and $b=2^{n-2}+1$. But then
$$
N+4D > b^m \geq b^2 =(2^{n-2}+1)^2 > 2^{n+1} > N+4D,
$$
a contradiction. If, however, we have $y_1=y_2$ and $x_3=n-2$, 
$(N+D)+(N+3D)=2(N+2D)$ implies that
$$
2^{n-1} +b^{y_1}=b^m,
$$
contradicting $m \geq 2$, unless $(b,m,y_1,n)=(3,2,0,4)$, which does not correspond to a $5$-term arithmetic progression in $S_{2,3}$.

We may therefore assume that $x_3=n$, so that, from $2(N+2D)=(N+D)+(N+3D)$, 
$$
2 \cdot b^{y_2}=2^{n-2}+b^{y_1}+b^{y_3}.
$$
Since $y_1 \geq y_2$, this implies that either
$$
2 \cdot b^{y_2} >  b \cdot b^{y_2},
$$
an immediate contradiction, or that $y_1=y_2$, whence
$$
b^{y_2}=2^{n-2}+b^{y_3},
$$
and, via Lemma \ref{tech-lem}, $y_3=0$ and either
$$
b=2^{n-2}+1, \; y_2=1,
$$
or
$$
(b,y_2,n) \in \{ (3,1,3), (3,2,5), (9,1,5) \}.
$$
The latter cases lead to no new arithmetic progressions while the former implies that
$$
N+4D > b^m \geq b^2 =(2^{n-2}+1)^2 > 2^{n+1}>N+4D,
$$
provided $n \geq 6$.

\subsection{The cases $x_1 \leq n-3$}

If we have $x_1 \leq n-3$, then 
$$
N+4D < 4(N+D) = 2^{x_1+2} + 4 \cdot b^{y_1} < \left( \frac{1}{2^{n-x_1-2}} + \frac{4}{b^{m-y_1}} \right) (N+4D)
$$
and hence $y_1=m-1$ and $b \leq 7$. Further, either $b \in \{ 3, 4 \}$, or $b=5$ and $x_1 \in \{ n-4, n-3 \}$, or $b \in \{ 6, 7 \}$ and $x_1=n-3$. If $N+3D$ is $2$-large, the identity $(N+D)+(N+3D)=2(N+2D)$  leads to the equation
$$
2 \cdot b^{y_2}=2^{x_1}+b^{m-1}+b^{y_3},
$$
an immediate contradiction modulo $3$ if $b \in \{ 4, 7 \}$, or if $y_2 < m-1$. We thus have
$$
b^{m-1} =2^{x_1}+b^{y_3},
$$
so that $y_3=0$ and 
$$
(b,m,x_1) \in \{ (3,2,1), (3,3,3), (5,2,2), (9,2,3) \}.
$$
None of these correspond to a new progression. It follows that 
$N+3D$ is $b$-large.
If $N+4D$ is $2$-large, then
$$
N=2(N+2D)-(N+4D) = 2 \cdot b^{y_2}-b^{y_4} > 1
$$
implies that $y_2 \geq y_4$ and that $y_2>0$.
Further, the identity $3(N+2D)=(N+4D)+2(N+D)$ implies
\begin{equation} \label{pig}
2^{n-1}+3 \cdot b^{y_2}=b^{y_4} +2^{x_1+1}+2 \cdot b^{m-1}.
\end{equation}
If $x_1=n-3$, it follows that
$$
2^{n-2}+3 \cdot b^{y_2}=b^{y_4}+2 \cdot b^{m-1},
$$
whence, modulo $3$, $b \in \{ 3, 5, 6 \}$ and, from $y_2>0$, necessarily $y_4=0$. The identity $3(N+3D)=2(N+4D)+(N+D)$ thus yields 
\begin{equation} \label{pig2}
3 \cdot 2^{x_3}+(3b-1) \cdot b^{m-1}=2 +9 \cdot 2^{n-3},
\end{equation}
a contradiction modulo $3$ for $b \in \{ 3, 6 \}$. We thus have $b=5$ and so, from (\ref{pig}), we have that $n \equiv 2 \mod{4}$, whence, considering (\ref{pig2}) modulo $5$, $x_3 \equiv 3 \mod{4}$. Since $(N+D)+(N+4D)=(N+2D)+(N+3D)$ implies the equation
$$
5 \cdot 2^{n-3}+1=5^{y_2}+2^{x_3}+4 \cdot 5^{m-1},
$$
we therefore have a contradiction modulo $5$.

If $x_1=n-3$ and both $N+3D$ and $N+4D$ are $b$-large, $x_4 \leq n-3$ implies that $D < 2^{n-3}$ and so
$$
N+4D=(N+D)+3D < 2^{n-3}+b^{m-1} + 3 \cdot 2^{n-3} <
\left( \frac{1}{2} + \frac{1}{b} \right) (N+4D),
$$
a contradiction for $b \geq 3$. Thus $x_4 \in \{ n-2, n-1 \}$.
If $x_4=n-2$, then, from $3(N+2D)=(N+4D)+2(N+D)$,
$$
2^n=b^{m}+2 \cdot b^{m-1}-3 \cdot b^{y_2}, 
$$
so that, modulo $3$, $b=5$ and $y_2=0$. But then, from $2(N+2D)=(N+D)+(N+3D)$,
$$
5 \cdot 2^{n-3} =6 \cdot 5^{m-1}-2,
$$
a contradiction. We thus have $x_4=n-1$ and so, from $(N+D)+(N+4D)=(N+2D)+(N+3D)$,
$$
b^{m-1}+2^{n-3}=b^{y_2}+2^{x_3}.
$$
If $x_3=n-2$, we therefore have
$$
b^{m-1}=b^{y_2}+2^{n-3}
$$
and so, from Lemma \ref{tech-lem} and $b \leq 7$, $y_2=0$ and 
$$
(b,m,n) \in \{ (3,2,4), (3,3,6), (5,2,5) \},
$$
none of which correspond to new progressions. We thus have $x_3=n-3$ and $y_2=m-1$, whereby $2(N+2D)=(N+D)+(N+3D)$ yields
$$
3 \cdot 2^{n-2} =(b-1) b^{m-1}
$$
and so $b=4$ and $n=2m$, which leads to $N=0$, a contradiction.. 

It remains, then, to consider the case where $x_1 \leq n-4$ and both $N+3D$ and $N+4D$ are $b$-large (so that $b \in \{ 3, 4, 5 \}$). As previously, we have $x_4 \in \{ n-2, n-1 \}$, so that
$$
2^{n-2}+b^m \leq 2^{x_4}+b^{m} =N+4D<4(N+D) \leq 2^{n-2}+ 4 \cdot b^{m-1}
$$
and so $b=3$. The identity $3(N+3D)=2(N+4d)+(N+D)$ thus implies
$$
2 \cdot 3^{m}=2^{x_4+1}+2^{x_1}-3 \cdot 2^{x_3},
$$
and so either $x_1=x_3=0$, or $\min \{ x_1, x_3 \}=1$. In either case, Proposition \ref{DT} and a short calculation completes the proof.

%---------------------------------------------------------------
\section{One large term} \label{Sec7}
%---------------------------------------------------------------

The last case to consider is when only the term $N+4D$ is large.
For this to happen, we require 
$$
N+3D \leq a^{n-1}+b^{m-1} < \left( \frac{1}{a} + \frac{1}{b} \right) (N+4D)
$$
and so $(a,b)=(2,3)$. More precisely, $N+3D=2^{n-1}+3^{m-1}$ and either $x_2=n-1$ or $y_2=m-1$, since otherwise we would have
$$
N+4D < 2(N+2D) \leq 2\left( 2^{n-2}+3^{m-2} \right) < \left( \frac{1}{2} + \frac{2}{9} \right) (N+4D).
$$
From $2(N+3D)=(N+2D)+(N+4D)$, we  have
$$
2^n+2 \cdot 3^{m-1} = 2^{x_2}+3^{y_2}+2^{x_4}+3^{y_4}
$$
and hence one of
$$
2 \cdot 3^{m-1} = 2^{n-1}+3^{y_2}+3^{y_4},
$$
$$
2^{n-1} = 3^{y_2}+2^{x_4}+3^{m-1},
$$
$$
3^{m-1} = 2^{x_2}+3^{y_4},
$$
or
$$
2^n = 2^{x_2}+2^{x_4}+2 \cdot 3^{m-1},
$$
depending on whether $N+4D$ is $2$-large or $3$-large. In each case, Proposition \ref{DT} and a little work reveal no new $5$-term arithmetic progressions in $S_{2,3}$. This completes the proof of Theorem \ref{thm-main}.

%--------------------------------------. 
\section{Six terms : the proof of Corollary \ref{cor-main}} \label{Sec8}
%--------------------------------------

If $S_{a,b}$ contains a $6$-term arithmetic progression, then, from Theorem \ref{thm-main}, we require $N+5D \in S_{a,b}$ for 
$$
(a,b,N,D)=(2,2^k+1,2^k+1,2^k)  \; \mbox{ or } \;  (3,4 \cdot 3^{k-1}+1,3^{k-1}+1,2 \cdot 3^{k-1}) \; \mbox{ for } \; k \in \mathbb{N},
$$
or
$$
\begin{array}{c}
(a,b,N,D) \in \left\{ (2,3,5,2),  (2,3,7,6), (2,3,9,8),  (2,3,17,24),  \right. \\
\left. (2,3,41,24), (2,5,5,8),  (2,9,17,24),  (2,9,41,24), (3,4,7,6) \right\}.
\end{array}
$$
A short check of the latter $9$ sporadic cases, reveals the $6$-term progressions corresponding to 
$$
(a,b,N,D) = (2,3,17,24) \; \mbox{ and } \; (2,9,17,24),
$$
and no $7$-term progressions. Suppose that
$(a,b,N,D)=(2,2^k+1,2^k+1,2^k)$, for $k$ a positive integer, and that 
$$
N+5D = 2^k+1 + 5 \cdot 2^k =3 \cdot 2^{k+1}+1= 2^{x_5} + (2^k+1)^{y_5},
$$
for $x_5$ and $y_5$ nonnegative integers. If $k=1$, then $x_5=y_5=2$, so that $(a,b,N,D) = (2,3,3,2)$. If $k \geq 2$, then 
$$
(2^k+1)^2 \geq 3 \cdot 2^{k+1}+1,
$$
whence necessarily $y_5 \in \{ 0, 1 \}$. Thus
$$
3 \cdot 2^{k+1} = 2^{x_5} \; \mbox{ or } \; 2^{k+2}+1= 2^{x_5},
$$
each an immediate contradiction.
Finally, suppose that $(a,b,N,D)=(3,4 \cdot 3^{k-1}+1,3^k+1,2 \cdot 3^k)$ for $k$ a positive integer, and that 
$$
N+5D = 3^{k-1}+1+ 10 \cdot 3^{k-1} = 11 \cdot 3^{k-1}+1=3^{x_5}+(4 \cdot 3^{k-1}+1)^{y_5},
$$
with $x_5$ and $y_5$ nonnegative integers. For each positive integer $k$, we have
$$
(4 \cdot 3^{k-1}+1)^2 > 11 \cdot 3^{k-1}+1
$$
and so $y_5 \in \{ 0, 1 \}$, corresponding to
$$
11 \cdot 3^{k-1}=3^{x_5} \; \mbox{ and } \;
7 \cdot 3^{k-1}=3^{x_5},
$$
respectively. These contradictions complete the proof of 
Corollary \ref{cor-main}.

%----------------------------------------------------------
\section{Infinite families} \label{Sec9}
%----------------------------------------------------------

While Theorem \ref{thm-main} provides a good understanding of $5$-term progressions in the sets $S_{a,b}$, it appears to be significantly more challenging to completely characterize all $3$ or $4$-term progressions. A  $3$-term arithmetic progression $N+iD = a^{x_i}+b^{y_i}$, $i \in \{ 0, 1, 2 \}$ in $S_{a,b}$ corresponds to an equation of the shape
\begin{equation} \label{fundy}
a^{x_0}+b^{y_0}+a^{x_2}+b^{y_2} = 2 \cdot a^{x_1}+ 2 \cdot b^{y_1},
\end{equation}
which, as a $6$-term $S$-unit equation (here, $S$ is the set of primes dividing $2ab$), has, typically,  at most finitely many solutions in integer exponents $x_i, y_i$, via the {\it Fundamental theorem of $S$-unit equations} (a consequence of Schmidt's Subspace Theorem) : 

\begin{thm}[Evertse \cite{Ev}] \label{S}
Let $S=\{p_1, p_2, \ldots, p_t \}$ be a finite set of primes, and $n$ a nonnegative integer. Then there are at most  finitely many integers $x_0, x_1, \ldots, x_n$, all of whose prime factors lie in $S$, and satisfying
\begin{equation} \label{sunit}
x_0+x_1+ \cdots + x_n=0,
\end{equation}
with $\gcd (x_0,x_1, \ldots, x_n)=1$ and without {\it vanishing subsums}, i.e. with $x_{i_1}+x_{i_2}+ \cdots + x_{i_k} \neq 0$, for every nontrivial, proper subset $\{ i_1, i_2, \ldots, i_k \}$ of $\{ 0, 1,  \ldots, n \}$.
\end{thm}
It follows that the question of  finiteness of $3$-term arithmetic progressions in $S_{a,b}$ is easily resolved unless either
\begin{enumerate}
\item the terms in (\ref{fundy}) have a common factor, or
\item equation (\ref{fundy}) has a vanishing subsum.
\end{enumerate}
These caveats, however, can lead to some surprisingly difficult arithmetic questions. By way of example, if, for a given integer $b>2$, we have even a single solution to the Diophantine equation
\begin{equation} \label{kruk}
1+b^{y_2} +2^{x_0} = 2 b^{y_1},
\end{equation}
in nonnegative integers $x_0, y_1$ and $y_2$, then $S_{2,b}$ has infinitely many $3$-term arithmetic progressions, given by
$$
N=2^{x_0}+1, \; N+D=2^{x_1}+b^{y_1}, \; N+2D=2^{x_1+1}+ b^{y_2},
$$
where $x_1$ can be taken to be any nonnegative integer. This corresponds to a vanishing subsum in (\ref{fundy}), specifically to the case where $a^{x_2}=2 \cdot a^{x_1}$. Presumably, equation (\ref{kruk}) has only  solutions when $b=2^k+1$ for some positive integer $k$, where we find
$$
(x_0,y_1,y_2) \in \{ (k,1,1), (k+1,1,0) \}
$$
and, additionally, in case $k=1$,
$$
(b,x_0,y_1,y_2) \in \{ (3,3,2,2), (3,4,2,0) \}.
$$
While finiteness for equation (\ref{kruk}) is a reasonably easy consequence of the abc-conjecture, it does not appear to follow unconditionally in an obvious way from, for example,  effective technology from Diophantine approximation, such as bounds for linear forms in logarithms. We should note, though, that, ineffectively, finiteness of the number of $3$-term arithmetic progressions in $S_{a,b}$ is almost immediate from Theorem \ref{S}, for any pair $(a,b)$ with $b>a>2$.

In the remainder of this section, we will discuss the various families of which we are aware for which $S_{a,b}$ is guaranteed to have infinitely many $3$-term arithmetic progressions, or any $4$-term arithmetic progressions.

\subsection{$3$-term arithmetic progressions in $S_{a,b}$}

If we take $a=2$ and $b=2^k+1$ for $k$ a positive integer, and set 
$N=2^k+1$ and $D=2^j$, where $j$ is a nonnegative integer,
then we have
$N+iD = a^{x_i}+b^{y_i}$,
with 
$$
(x_0,y_0) = (k,0), \; (x_1,y_1)=(j,1) \; \mbox{ and } \; (x_2,y_2) =
(j+1,1),
$$
so that $S_{2,2^k+1}$ contains infinitely many $3$-term arithmetic progressions. Similarly, we find a second family by taking $N=2^{k+1}+1$ and $D=2^j-2^k$ for $j \geq k+1$. In this case
we have
$N+iD = a^{x_i}+b^{y_i}$,
with 
$$
(x_0,y_0) = (k+1,0), \; (x_1,y_1)=(j,1) \; \mbox{ and } \; (x_2,y_2) =
(j+1,0).
$$

Suppose now that $a$ and $b$ are {\it multiplicatively dependent}, say with $a^c=b^d$ for $c$ and $d$  positive integers.
Then setting $N=a^{kc}+b^{kd}$ and $D=a^{(k+j)c}-a^{kc}$ gives $N+iD = a^{x_i}+b^{y_i}$,
with 
$$
(x_0,y_0) = (kc,kd), \; (x_1,y_1)=((k+j)c,kd) \; \mbox{ and } \; (x_2,y_2) =((k+j)c,(k+j)d),
$$
for any nonnegative integer $k$ and positive integer $j$, and hence infinitely many $3$-term arithmetic progressions in $S_{a,b}$.

\subsection{$4$-term arithmetic progressions in $S_{a,b}$}

Continuing with the theme of multiplicatively dependent $a$ and $b$, suppose that
we have 
\begin{equation} \label{spek}
(a,b) = (2^d, 2^c ) \mbox{ with } d < c \mbox{ positive integers and } \gcd (d, c)=1.
\end{equation}
Then $a$ and $b$ are multiplicatively dependent (with $a^c=b^d$) and, 
if $k$ and $j$ are positive integers with $dk-cj=1$,  setting
$$
N = 2^{kd} + 2^{jc} \; \mbox{ and }  \; D=2^{jc + m cd }-2^{j c}, \; \mbox{ for } m \geq 1,
$$
we find a $4$-term arithmetic progression in $S_{a,b}$, given by the terms
$$
2^{k d} + 2^{j c}, \; 2^{j d + m cd } + 2^{k d}, \; 2^{k d + m cd } + 2^{jc} \; \mbox{ and }  \; 2^{k d + m cd } +2^{j c + m cd }.
$$
Similarly, if $k$ and $j$ are positive integers with $d k-c j=-1$, and we set
$$
N = 2^{k d} + 2^{jc}, \; \; D=2^{kd + m cd }-2^{k d}, \; \mbox{ for } m \geq 1,
$$
then the terms
$$
2^{k d} + 2^{j c}, \; 2^{k d + m cd } + 2^{j c}, \; 2^{j c + m cd } + 2^{k d} \; \mbox{ and } \; 2^{k d + m cd } +2^{j c + m cd }
$$
form a $4$-term arithmetic progression in $S_{a,b}$. It follows that there are infinitely many $4$-term arithmetic progressions in $S_{a,b}$.

\subsection{Families of $(a,b)$ for which $S_{a,b}$ contains a $4$-term arithmetic progression}

While the pairs $(a,b)$ given in (\ref{spek}) are likely the only ones for which $S_{a,b}$ has infinitely many $4$-term arithmetic progressions, 
there are a number of other infinite families of pairs $(a,b)$ (each with $\gcd (a,b)=1$) for which the set of $4$-term arithmetic progressions in $S_{a,b}$ is nonempty, even excluding those coming from the $5$-term families in Theorem \ref{thm-main}. These correspond to the following terms forming $4$-term arithmetic progressions in $S_{a,b}$ :
\begin{equation} \label{prog1}
a^0+b^0, \; a^1+b^0, \; a^0+b^1 \;   \mbox{ and }  \;  a^1+b^1,
\end{equation}
\begin{equation} \label{prog2}
a^0+b^0, \; a^0+b^1, \; a^t+b^0  \;   \mbox{ and }  \;  a^t+b^1, \; t \geq 2,
\end{equation}
\begin{equation} \label{prog3}
a^{\delta_1}+b^{\delta_2}, \; a^2+b^{\delta_2}, \; a^{\delta_1}+b^2 \;   \mbox{ and }  \;  a^2+b^2, \; \delta_1, \delta_2 \in \{ 0, 1 \},
\end{equation}
\begin{equation} \label{prog4}
a^{t+1}+b^0, \; a^{t+1}+b^2, \; a^{2t+1}+b^0 \;  \mbox{ and }  \; a^{2t+1}+b^2,  \; t \geq 2,
\end{equation}
\begin{equation} \label{prog5}
a^{t+3}+b^0, \; a^{t+3}+b^2, \; a^{2t+3}+b^0 \;  \mbox{ and }  \; a^{2t+3}+b^2,  \; t \geq 2,
\end{equation}
\begin{equation} \label{prog6}
a^{t+1}+b^1, \; a^{t}+b^2, \; a^{2t}+b^2 \;  \mbox{ and }  \; a^{3t}+b^1,  \; t \geq 1
\end{equation}
and
\begin{equation} \label{prog7}
a^{s}+b^1, \; a^{s+1}+b^1, \; a^{t}+b^0 \;  \mbox{ and }  \; a^{s+2}+b^1,  \; 1 \leq s \leq t-2.
\end{equation}
In case (\ref{prog1}), it suffices to take
$$
(a,b,N,D) = (n,2n-1,2,n-1),  \mbox{ for $n \geq 2$ an integer.}
$$
For the terms in (\ref{prog2}) to form a progression, we require
$$
(a,b,N,D) =( 2k+1,(a^t+1)/2,2,(a^t-1)/2),  \mbox{ for $k \geq 1$ and $t \geq 2$  integers.}
$$
In the case of (\ref{prog3}), 
necessarily
$$
\begin{array}{l}
(a,b,N,D) = (a,b,a^{\delta_1}+b^{\delta_2},a^2-a^{\delta_1}),  \mbox{ for integers } \\
b>a>1 \mbox{ with $b^2-b^{\delta_2}=2a^2-2a^{\delta_1}, \; \delta_1, \delta_2 \in \{ 0, 1 \}$}. \\
\end{array}
$$
Note that this last equation may be  readily shown to have infinitely many solutions in integers, for any of the four choices of $(\delta_1,\delta_2)$. For (\ref{prog4}) and (\ref{prog5}), we have
$$
(a,b,N,D)=(8,2^{3t+1}-1,8^{t+1}+1, b^2-1)
$$
and
$$
(a,b,N,D)=(2,2^{t+1}-1,2^{t+3}+1, b^2-1), 
$$
respectively, where $t$ is a positive integer, while, for (\ref{prog6}), we have
$$
(a,b,N,D)=(3,3^t+1,4 \cdot 3^t+1,3^{2t}-3^t), \; t \geq 1.
$$
Finally, in case (\ref{prog7}), 
$$
(a,b,N,D)=(2,2^t-3 \cdot 2^s+1,2^t-2^{s+1}+1,2^s), \mbox{ with } 1 \leq s \leq t-2.
$$

 There are possibly more such families to discover, as well as sporadic examples of pairs $(a,b)$ for which  $S_{a,b}$ contains $4$-term arithmetic progressions $N+iD$, $i \in \{ 0, 1, 2, 3 \}$. We know of few examples, beyond those previously discussed, for which $N > a+b$. An amusing one occurs for $(a,b)=(22,78)$, where 
$$
22+78^2, \; 22^4+78^2, \; 22+78^3 \; \mbox{ and } \; 22^4+78^3
$$
forms an arithmetic progression.

%----------------------------------
\section{Concluding remarks}
%----------------------------------

Besides the problem of characterizing $3$ or $4$-term arithmetic progressions in the $S_{a,b}$, it would be interesting to understand, for example, progressions in  sumsets of the shape
$$
S_{a,b,c} = \{ n \in \mathbb{N} \; : \; n = a^x+b^y+c^z, \; \mbox{ for } x, y, z \in \mathbb{Z}, \; x, y, z \geq 0 \},
$$
where, say, $c \geq b \geq a \geq 2$. Again, the machinery of $S$-units allows one, in certain cases, to make qualitative statements, but explicit results are apparently rather difficult to obtain.

%----------------------------------

\end{document}